\documentclass{amsart}

\usepackage[T1]{fontenc}
\usepackage[utf8]{inputenc}
\usepackage{lmodern}
\usepackage[a4paper,margin=28mm]{geometry}
\usepackage{amsmath,amssymb,amsthm,mathtools}
\usepackage{enumitem}
\usepackage{microtype}
\usepackage[hidelinks]{hyperref}
\usepackage[numbers,sort&compress]{natbib}
\numberwithin{equation}{section}
\newtheorem{theorem}{Theorem}[section]
\newtheorem{proposition}[theorem]{Proposition}
\newtheorem{lemma}[theorem]{Lemma}
\newtheorem{corollary}[theorem]{Corollary}
\theoremstyle{definition}
\newtheorem{definition}[theorem]{Definition}
\theoremstyle{remark}
\newtheorem{remark}[theorem]{Remark}
\newcommand{\FF}{\mathbb F_2}
\newcommand{\ZZ}{\mathbb Z_2}
\newcommand{\QQ}{\mathbb Q_2}
\newcommand{\CC}{\mathcal C}
\newcommand{\gr}{\operatorname{gr}}
\newcommand{\LL}{\widehat{\mathfrak L}}
\newcommand{\Prim}{\operatorname{Prim}}
\newcommand{\ad}{\operatorname{ad}}
\newcommand{\Ad}{\operatorname{Ad}}
\newcommand{\Sch}{\operatorname{Sch}}
\newcommand{\im}{\operatorname{im}}
\newcommand{\id}{\operatorname{id}}
\newcommand{\eps}{\epsilon_{\varnothing}}
\newcommand{\shuffle}{\mathbin{\amalg}}
\newcommand{\GA}{\mathbf A}
\newcommand{\GB}{\mathbf B}
\newcommand{\GD}{\mathbf D}
\newcommand{\GOne}{\mathcal G_1}
\newcommand{\GTwo}{\mathcal G_2}
\newcommand{\Ima}{\mathcal I}
\newcommand{\cR}{\mathcal R}
\newcommand{\cS}{\mathcal S}
\newcommand{\cQ}{\mathcal Q}
\newcommand{\cE}{\mathcal E}
\newcommand{\rLead}{\mathcal R_{\mathrm{lead}}}
\newcommand{\Slin}{\mathcal S_{\mathrm{lin}}}
\DeclareMathOperator{\GT}{GT}
\DeclareMathOperator{\Lie}{Lie}
\DeclareMathOperator{\rev}{rev}

\begin{document}

\title{A remark on the AI proof of the Deligne--Drinfeld conjecture}

\author{Nikita Markarian}

\begin{abstract}
We give a geometric interpretation of the projection and
reconstruction underlying the characteristic-two proof of the
Deligne--Drinfeld conjecture in \cite{DD}.
Using the double covering \(z\mapsto z^2\) and filling a puncture,
we construct a morphism on normalized level-two associators with
an explicit inverse on its image.
We compute the restricted tangent reconstruction in the pro-\(2\)
realization and identify the leading coefficient maps with those
of \cite{DD}, making its implicit reconstruction explicit.
\end{abstract}

\email{nikita.markarian@gmail.com}

\date{}

\address{UMR 7501, Université de Strasbourg,
7 rue René Descartes,
67084 Strasbourg Cedex, France}
\maketitle

\section*{Introduction}

The characteristic-two argument in the proof of the
Deligne--Drinfeld conjecture \cite{DD} uses an injective projection
of infinitesimal associators. Willwacher gives a simplified
exposition of the proof in \cite{Willwacher}. In this note, we interpret this
projection geometrically through the double covering \(z\mapsto z^2\)
and the operation of filling a puncture. We also give an explicit
formula for the reconstruction implicit in the injectivity proof
of \cite{DD}.

Filling the puncture \(0\) in
\(\mathbb P^1\setminus\{0,1,-1,\infty\}\) defines a morphism
from the scheme of normalized level-two associators to a scheme of
group-like series. We construct its inverse on the image, which
we call the Hirose--Sato map. Enriquez's isomorphism relates the
level-one and level-two pro-\(2\) deformation groups through the
double covering. In this realization, we compute the restricted
tangent reconstruction and show that the appropriate leading parts
of the forward and reconstruction coefficient formulas recover
the maps of \cite{DD}.

Section~\ref{sec:preliminaries} fixes the conventions and recalls
associators and the filtrations used below.
Sections~\ref{sec:enriquez-map}--\ref{sec:reconstruction} give the
normalized covering map and the rational and Betti inverse formulas.
Section~\ref{sec:restricted-tangent} proves strictness of the
cyclotomic projection and computes its restricted tangent inverse,
using the integral coefficient calculation of
Appendix~\ref{app:coefficient-proof}.

Section~\ref{sec:dd-comparison} compares the leading coefficient
maps at fixed covering weight with those of \cite{DD}, separately
from the intrinsic upstairs strictness of
Section~\ref{sec:restricted-tangent}. Ordinary duality and Lie
membership alone bound the projected image; transferring this
bound to the source requires the pentagon-dependent statement
that every nonzero leading source component survives projection.
We identify its reconstruction
operator with the leading part of Willwacher's pentagon formula.
Appendix~\ref{app:dd-reconstruction} derives the same operator
from the deletion tests and special identity of \cite{DD}.

Appendix~\ref{app:three-relations} studies the larger solution set
on the normalized covering graph defined by the three Betti
relations. It gives successive finite lifting tests for its
realized initial coefficients and a first example of an obstruction.
This supplementary analysis does not enter the comparison theorem.

\section{Conventions and associators}
\label{sec:preliminaries}

We use the field \(\FF\), the ring \(\ZZ\) of \(2\)-adic
integers, and its fraction field \(\QQ\).
Rational associator identities are considered
over commutative \(\mathbb Q\)-algebras; comparison with pro-\(2\)
groups uses \(\QQ\). All deformation groups have scalar parameter
\(1\). The associator parameter \(\mu\) is kept explicit.

Let \(\Bbbk\) denote the coefficient ring. For a finite alphabet \(S\), let \(\Bbbk\langle S\rangle\) be the
free associative algebra and \(\Lie_{\Bbbk}(S)\) its ordinary free Lie
subalgebra, with bracket \([u,v]=uv-vu\). Lie brackets are
alternating, including in characteristic two. Completions are in
positive weight; letters have weight one unless otherwise stated.
Write \(\LL_{\Bbbk}(S)\) for the completed ordinary free Lie algebra.
The tensor coproduct makes the letters primitive. Over a
\(\mathbb Q\)-algebra, the group-like series with constant term
one form \(\exp\LL_{\Bbbk}(S)\). Over \(\FF\), the primitives form
the completed free restricted Lie algebra, with restricted square
\(u^{[2]}=u^2\) \cite[Sections~4--5]{Passman}. A restricted Lie algebra over a field of characteristic two is a Lie
algebra with a square operation satisfying
\begin{equation}
 (\lambda u)^{[2]}=\lambda^2u^{[2]},\qquad
 \ad_{u^{[2]}}=\ad_u^2,\qquad
 (u+v)^{[2]}=u^{[2]}+v^{[2]}+[u,v].
 \label{eq:restricted-identities}
\end{equation}
\begin{remark}
Ordinary Lie membership is stronger than primitivity over \(\FF\).
We use restricted primitives for the group-associated graded, and
ordinary Lie polynomials for the reduced rational solutions of \cite{DD}.
\end{remark}

For a word \(w\), let \(c_w(\Omega)\) be its coefficient in a series
\(\Omega\), and set \(\lambda_\Omega(w)=c_w(\Omega)\). With the tensor
coproduct just specified, group-like series correspond to unital
characters of the shuffle algebra. A primitive series has zero
constant coefficient and its coefficient functional vanishes on
shuffle products of two nonempty words. The empty-word coefficient
is denoted by \(\eps\). Write \(\rev\) for word reversal,
and \(L_u,R_u\) for left and right multiplication by \(u\).
Operator products act rightmost first. The antipode sends a word
\(w\) to \((-1)^{|w|}\rev(w)\) and a primitive \(p\) to
\(-p\); hence reversal fixes every primitive over \(\FF\),
including restricted powers.

Group commutators are \([g,h]=ghg^{-1}h^{-1}\). Magnus
coordinates have the form \(g-1\); their Hopf coproduct is
\(\Delta(g-1)=(g-1)\otimes1+1\otimes(g-1)+(g-1)\otimes(g-1)\).
The primitive-letter tensor coproduct describes the associated
graded used below, rather than these ungraded coordinates.
For a free pro-\(2\) group \(G\), let
\begin{equation}
 D_nG=\{g\in G:g-1\in I_G^n\},
 \qquad I_G=\ker\bigl(\FF[[G]]\to\FF\bigr).
 \label{eq:zassenhaus}
\end{equation}
Here the group algebra and \(D_\bullet G\) use the ordinary
free-group law. This is that group's intrinsic Zassenhaus filtration.

For background on the Zassenhaus filtration and the restricted
Lie structure on its associated graded, see
\cite[Proposition~6.4 and Section~6.5]{DarneSuciu}.
Functoriality follows from the preservation of augmentation
ideals and their powers under the induced homomorphisms of
completed group algebras.

\begin{definition}[Restricted central filtrations]
A restricted central filtration on a group \(G\) is a decreasing
sequence of subgroups \(F_1G=G,F_2G,\ldots\) such that
\[
 [F_mG,F_nG]\subset F_{m+n}G,\qquad
 g^2\in F_{2n}G\quad(g\in F_nG).
\]
It is separated if \(\bigcap_nF_nG=\{1\}\). Its associated
graded
\[
 \gr_FG=\bigoplus_{n\ge1}F_nG/F_{n+1}G
\]
is a restricted Lie algebra over \(\FF\): multiplication gives
addition on each quotient, while commutators and squaring induce
the bracket and restricted square. A filtered homomorphism
\(\theta:G\to G'\) induces \(d\theta:=\gr_F\theta\).
It is strict if
\(\theta(F_nG)=\theta(G)\cap F_nG'\) for every \(n\).
\end{definition}

\begin{remark}[Restricted tangent spaces]
We call \(\gr_FG\) the restricted tangent space for the chosen
filtration; this is not the space of solutions over dual numbers.
We use the intrinsic upstairs Zassenhaus filtration in
Section~\ref{sec:restricted-tangent}, and the covering-weight
filtration \(W\) in Section~\ref{sec:weighted-enriquez}.
The deformation groups have composition law \(\star\), which
induces their restricted Lie operations. Brackets in coefficient
formulas are the usual tensor-algebra brackets.
\end{remark}

The infinitesimal braid Lie algebra \(\mathfrak t_n\) has generators
\(t_{ij}=t_{ji}\), \(i\ne j\), and relations
\[
 [t_{ij},t_{kl}]=0\quad(i,j,k,l\text{ distinct}),\qquad
 [t_{ij},t_{ik}+t_{jk}]=0\quad(i,j,k\text{ distinct}).
\]
\begin{definition}[Level-one associators]
For a commutative \(\mathbb Q\)-algebra \(\Bbbk\) and \(\mu\in\Bbbk\),
let \(\CC_{1,\mu}(\Bbbk)\) be the set of ordinary associators with
parameter \(\mu\). An element is a group-like series \(\Phi(u,v)\), with zero linear
coefficients, satisfying
\begin{align}
 \Phi(t_{12},t_{23}+t_{24})\Phi(t_{13}+t_{23},t_{34})
 &=\Phi(t_{23},t_{34})\Phi(t_{12}+t_{13},t_{24}+t_{34})
   \Phi(t_{12},t_{23}),\label{eq:ordinary-pentagon}\\
 \Phi(u,v)\Phi(v,u)&=1,\label{eq:ordinary-duality}\\
 e^{\mu u/2}\Phi(w,u)e^{\mu w/2}\Phi(v,w)e^{\mu v/2}\Phi(u,v)
 &=1,\qquad w=-u-v.\label{eq:ordinary-hexagon}
\end{align}
The first identity is in the completed enveloping algebra of
\(\mathfrak t_4\).
\end{definition}

Linearizing at \(\Phi=1\), with \(\mu=0\) fixed, gives
the ordinary infinitesimal associator equations. Their homogeneous
rational solution space in weight \(n\) is denoted by
\(\mathfrak{grt}_{1,n}\). In Section~\ref{sec:dd-comparison},
we use the lattice and its reduction
\[
 \Lambda_n=\mathfrak{grt}_{1,n}\cap
 \Lie_{\mathbb Z_{(2)}}(x,y)_n,\qquad
 L_n=\Lambda_n/2\Lambda_n,
 \qquad \mathbb Z_{(2)}=\{a/b\in\mathbb Q:b\text{ odd}\}.
\]
The intersection is taken in \(\Lie_{\mathbb Q}(x,y)_n\).
The lattice is saturated. Indeed, if \(2v\in\Lambda_n\) for an
integral \(v\), then \(v\) belongs to the same rational solution
space. Hence \(L_n\) embeds in \(\Lie_{\FF}(x,y)_n\).
Its elements satisfy the reduced equations and, for \(n>2\),
the reduced special (cabling) identity of
\cite[Lemma~2.4 and Proposition~2.5]{DD}.

Consider the arrangement in \(\mathbb C^3\) consisting of
\(z_i=0\), \(z_i=z_j\), and \(z_i=-z_j\).
Its integral holonomy Lie algebra has generators
\(a_i,p_{ij},m_{ij}\) for these hyperplanes. For every codimension-two
intersection \(L\) and every hyperplane \(H\supset L\), its
defining relation is
\[
 [t_H,\ \sum_{H'\supset L}t_{H'}]=0.
\]
We use its degree completion over \(\Bbbk\) for the mixed pentagon.
The notation is \(a_i=t^{1,i+1}\),
\(p_{ij}=t(0)^{i+1,j+1}\), \(m_{ij}=t(1)^{i+1,j+1}\).

\begin{lemma}[The arrangement augmentation algebra]
\label{lem:arrangement-graded-algebra}
Let \(\mathfrak h_{\mathbb Z}\) be the integral holonomy Lie
algebra of this arrangement, \(P_{\mathrm{arr}}=P(2,3)\) its
pure monomial braid group, and \(I\) the augmentation ideal of
\(\FF[[\widehat P_{\mathrm{arr}}^{(2)}]]\). For a small positive
loop \(g_H\) around \(H\), the assignment
\(t_H\mapsto (g_H-1)\bmod I^2\) induces an isomorphism of
graded Hopf algebras
\begin{equation}
 U(\mathfrak h_{\mathbb Z})\otimes\FF
 \xrightarrow{\ \sim\ }
 \gr_I\FF[[\widehat P_{\mathrm{arr}}^{(2)}]].
 \label{eq:arrangement-graded-algebra}
\end{equation}
\end{lemma}

\begin{proof}
The group is an iterated almost-direct product of free groups
of ranks \(1,3,5\). Its integral lower-central-series Lie algebra
is \(\mathfrak h_{\mathbb Z}\), and as a graded abelian group
is the sum of free Lie algebras on \(1,3,5\) generators
\cite[Theorem~2.3.6 and Corollary~2.3.9]{Cohen}.
In particular the holonomy relations hold for the initial loop classes
in the augmentation-graded algebra. They define a surjection from
the left-hand side of \eqref{eq:arrangement-graded-algebra}, since
these classes generate the latter algebra in degree one.

For the discrete group, \cite[Theorem~6.13]{DarneSuciu} applies:
the action at each stage is trivial on integral abelianization,
hence on its reduction modulo two. Thus its Zassenhaus restricted
Lie algebra is an iterated semidirect product of the corresponding
free restricted Lie algebras. The restricted enveloping-algebra theorem recalled in
\cite[Section~6.5]{DarneSuciu}, together with restricted PBW,
gives the augmentation-graded algebra the Hilbert series
\[
 \frac{1}{(1-t)(1-3t)(1-5t)}.
\]
Ordinary PBW and Cohen's integral decomposition give the same
series for \(U(\mathfrak h_{\mathbb Z})\otimes\FF\).

Passing to pro-\(2\) completion does not change the finite
augmentation quotients. Indeed, for a finitely generated group
\(G\), the algebra \(\FF[G]/I_G^n\) is finite-dimensional,
and \(1+I_G/I_G^n\) is a finite \(2\)-group. The natural
map from \(G\) into this group therefore extends to its pro-\(2\)
completion. Conversely, the augmentation ideal of the group
algebra of any finite \(2\)-group is nilpotent. These two facts
identify the augmentation completion of \(\FF[G]\) with
\(\FF[[\widehat G^{(2)}]]\) and identify their finite
augmentation quotients. The surjection above consequently has
equal finite dimensions in every degree and is an isomorphism.
It respects coproducts because the degree-one loop classes
are primitive.
\end{proof}

\begin{definition}[Level-two associators]
\label{def:cyclotomic-pair}
Put \(a_\infty=-a-b-d\), the residue at infinity.
We denote by \(\CC_{2,\mu}(\Bbbk)\) the set of normalized level-two
cyclotomic associators with parameter \(\mu\). Its elements are
pairs
\[
 (\Phi,\Psi)\in\exp\LL_{\Bbbk}(u,v)\times\exp\LL_{\Bbbk}(a,b,d),
\]
where \(\Phi\in\CC_{1,\mu}(\Bbbk)\), \(c_b(\Psi)=0\),
and \(\Psi\) satisfies the mixed pentagon
\begin{equation}
 \Psi^{1,2,34}\Psi^{12,3,4}
 =\Phi^{2,3,4}\Psi^{1,23,4}\Psi^{1,2,3}
 \label{eq:mixed-pentagon}
\end{equation}
and the octagon
\begin{align}
 &\Psi(a,d,b)^{-1}e^{\mu d/2}\Psi(a_\infty,d,b)e^{\mu a_\infty/2}
 \notag\\[-1mm]
 &\hspace{12mm}\cdot
 \Psi(a_\infty,b,d)^{-1}e^{\mu b/2}\Psi(a,b,d)e^{\mu a/2}=1.
 \label{eq:octagon}
\end{align}
\end{definition}

The superscripts in \eqref{eq:mixed-pentagon} are the standard
insertions into this completed holonomy algebra. The four insertions, in their order in
\eqref{eq:mixed-pentagon}, are
\begin{align*}
 \Psi^{1,2,34}&=\Psi(a_1,p_{12}+p_{13},m_{12}+m_{13}),\\
 \Psi^{12,3,4}&=\Psi(a_2+p_{12}+m_{12},p_{23},m_{23}),\\
 \Psi^{1,23,4}&=\Psi(a_1+a_2+p_{12}+m_{12},
                 p_{13}+p_{23},m_{13}+m_{23}),\\
 \Psi^{1,2,3}&=\Psi(a_1,p_{12},m_{12}),
 \qquad\Phi^{2,3,4}=\Phi(p_{12},p_{23}).
\end{align*}
These are the conventions of \cite[Section~1]{Furusho}.

\begin{proposition}[Specialization and duality]
\label{prop:boundary-duality}
We shall use two consequences of the normalized mixed pentagon:
\begin{equation}
 c_a(\Psi)=0,\qquad \Psi(a,b,0)=\Phi(a,b),
 \label{eq:ordinary-specialization}
\end{equation}
and Broadhurst duality
\begin{equation}
 \Psi(a,b,d)=e^{-\rho b}
 \Psi(b,a,-a-b-d)^{-1}e^{-\rho a},
 \qquad \rho=c_d(\Psi).
 \label{eq:broadhurst}
\end{equation}
\end{proposition}

\begin{proof}
In the mixed pentagon, retain only \(a_1=t\). This gives
\(\Psi(t,0,0)=\Psi(t,0,0)^2\), hence \(\Psi(t,0,0)=1\).
Normalization gives \(\Psi(0,t,0)=1\) as well. Now kill all
\(a_i,m_{ij}\) and put
\((p_{12},p_{23},p_{13})=(u,v,-u-v)\).
This substitution respects the holonomy relations since
\(p_{12}+p_{13}+p_{23}=0\). The mixed pentagon becomes
\[
 \Psi(u,v,0)=\Phi(u,v)\Psi(u,-u,0).
\]
The last factor is \(1\), since both linear coefficients vanish
and all Lie brackets vanish on a single generator. This proves
\eqref{eq:ordinary-specialization}. The pair now has the form
\((\Psi(a,b,0),\Psi)\) in Hirose's definition of
\(\mathrm{PENT}(2)\). His Proposition~24 gives
\eqref{eq:broadhurst} over every commutative \(\mathbb Q\)-algebra
\cite[Definition~11 and Proposition~24]{Hirose}.
\end{proof}

The functors \(\CC_{1,\mu}\) and \(\CC_{2,\mu}\) are defined
by coefficient equations in completed group-like series. We regard
them as affine schemes, without a finite-type assumption, over the
chosen parameter base \(\mathbb Q[\mu]\), or after fixing \(\mu\).
The target schemes below are understood over the same base. Write
\(\mathcal U(S)(\Bbbk)=\exp\LL_{\Bbbk}(S)\) for the corresponding
scheme of group-like series.

\begin{remark}[Coordinates and realizations]
The letters \(u,v\) and \(a,b,d\) denote residue variables
at levels one and two. Upstairs pro-\(2\) generators are bold;
their initial forms will be specified with the filtration.
The capitals \(A,B,C\) denote the covering coordinates of \cite{DD}.
\end{remark}

Let \(F=\widehat F^{(2)}(X,Y)\) be the free pro-\(2\) group.
\begin{definition}
The group \(\GOne=\GT^{(2)}_1\) consists of the parameter-one
pro-\(2\) Grothendieck--Teichm\"uller elements \(f\in F\).
Explicitly, \(f\) has trivial abelianization and satisfies
\begin{align}
 f(Y,X)f(X,Y)&=1,\label{eq:group-duality}\\
 f(Z,X)f(Y,Z)f(X,Y)&=1,
       \qquad Z=(XY)^{-1},\label{eq:group-threecycle}\\
 f(x_{12},x_{23}x_{24})f(x_{13}x_{23},x_{34})
 &=f(x_{23},x_{34})f(x_{12}x_{13},x_{24}x_{34})f(x_{12},x_{23}).
 \label{eq:group-pentagon}
\end{align}
The last equality is in \(\widehat P_4^{(2)}\), the pro-\(2\)
completion of the pure braid group. Here
\(x_{ij}\) are the Artin generators. For \(v\in F\), write
\(\theta_v(X)=vXv^{-1}\) and \(\theta_v(Y)=Y\). The group law
is \(f_1\star f_2=\theta_{f_2}(f_1)f_2\).
\end{definition}

\section{Enriquez's map}
\label{sec:enriquez-map}

The double covering \(z\mapsto z^2\) relates the level-one and
level-two pro-\(2\) fundamental groups. Enriquez's level-change
isomorphism identifies their normalized deformation groups. We fix
the covering generators and write this isomorphism explicitly.

In the free pro-\(2\) group \(F\) of Section~\ref{sec:preliminaries}, put
\[
 \Gamma=\ker\bigl(F\longrightarrow\mathbb Z/2\mathbb Z\bigr),
 \qquad X\longmapsto1,\quad Y\longmapsto0.
\]
It is freely generated, as a pro-\(2\) group, by
\begin{equation}
 \GA=X^2,\qquad \GB=Y,\qquad \GD=XYX^{-1}.
 \label{eq:schreier-generators}
\end{equation}
Write \(j:\Gamma\hookrightarrow F\) for the inclusion and
\(\Sch:j(\Gamma)\to\Gamma\) for its inverse, expressed in these
generators.
The exponent homomorphisms \(e_A,e_B,e_D:\Gamma\to\ZZ\)
refer throughout to the basis \((\GA,\GB,\GD)\).

Geometrically, \(j\) is induced by the covering
\[
 \mathbb P^1\setminus\{0,1,-1,\infty\}
 \longrightarrow\mathbb P^1\setminus\{0,1,\infty\},
 \qquad z\longmapsto z^2.
\]
We choose basepoints and paths so that the loops around \(0,1,-1\) are
\eqref{eq:schreier-generators}.

We use Enriquez's parameter-one level-two deformation group, with its
extra module parameter removed by the normalization \(e_B(h)=0\).
The following description also fixes the normalization of its
identification with \(\GOne\).

\begin{proposition}[Normalized Enriquez isomorphism]
\label{prop:enriquez}
For \(f\in\GOne\), set
\begin{equation}
 u_f=\Sch(f),\qquad s_f=-e_B(u_f),\qquad
 h_f=\GB^{s_f}u_f.
 \label{eq:enriquez-formula}
\end{equation}
Then
\begin{equation}
 \cE:\GOne\xrightarrow{\ \sim\ }\GTwo,
 \qquad f\longmapsto(f,h_f),
 \label{eq:enriquez}
\end{equation}
is an isomorphism onto the normalized level-two deformation group.
Its inverse forgets \(h\). Every \((f,h)\in\GTwo\) satisfies
\begin{equation}
 e_A(h)=e_B(h)=0,\qquad
 j(h)=Y^s f,\qquad s=e_D(h).
 \label{eq:covering-identity}
\end{equation}
\end{proposition}

\begin{proof}
Enriquez's level-change map removes the full \(2\)-part of
\(N\) in the pro-\(2\) category, so it identifies level two
with level one \cite[Section~6.3]{Enriquez}.
At level one let \(g\) be the module component and write
\(g(1,Y)=Y^s\). Forgetting the first strand in the mixed pentagon
and using the centrality of \(x_{12}x_{13}x_{23}\) gives
\(g(X,Y)=Y^s f(X,Y)\), as in the proof of
\cite[Proposition~5.4]{Enriquez}. Conversely, substituting
\(Y^s f\) into the mixed pentagon multiplies each side of the
ordinary pentagon on the left by \((x_{23}x_{24}x_{34})^s\)
in strand labels \(1,2,3,4\). Here the powers of \(x_{34}\)
and \(x_{23}\) commute past the adjacent factors by the disjoint
strand relations and the central full twists on triples.
The octagon follows from
ordinary duality and the three-cycle by the calculation in
Section~\ref{sec:reconstruction}. These identities hold for
\(s\in\ZZ\) by continuity. Since \(f\) has zero
\(X\)-exponent, it lies in
\(j(\Gamma)\). Its lift has \(e_A(u_f)=0\) and
\(e_B(u_f)+e_D(u_f)=0\). Multiplication by \(\GB^{s_f}\)
therefore gives exactly \eqref{eq:covering-identity} and the required
normalization. Conversely, that normalization uniquely determines
\(s\) from \(f\).

For the group assertion, the parameter-one action on
\(\Gamma^{\mathrm{ab}}\) is the identity. Hence \(e_B\) is
additive for the cyclotomic composition law. Its zero set is a
subgroup, and the normalized bijection is a group isomorphism.
\end{proof}

\begin{remark}[Normalization]
The exponent \(e_B\) always refers to the puncture basis
\((\GA,\GB,\GD)\). In the alternative basis
\((\GA,\GB,\mathbf C)\), where \(\mathbf C=\GD\GB^{-1}\),
one has \(e_B(\mathbf C)=-1\). Thus the normalization in \eqref{eq:enriquez-formula} uses
the puncture basis, even when covering coordinates are used later.
\end{remark}

\section{Cyclotomic projection and the Hirose--Sato map}
\label{sec:forward-maps}
\label{sec:reconstruction}

Filling the puncture \(0\) upstairs defines the cyclotomic
projection. We describe it on associator schemes and on pro-\(2\)
deformation groups, then construct its inverse on the image.

\begin{definition}[Cyclotomic projection]
\label{def:cyclotomic-projection}
The morphism of associator schemes considered here is
\begin{equation}
 q_{\mathrm{rat}}:\CC_{2,\mu}\longrightarrow\mathcal U(b,d),
 \qquad (\Phi,\Psi)\longmapsto \Psi(0,b,d).
 \label{eq:cyclotomic-projection}
\end{equation}
Write \(\Ima_{\mathrm{rat}}(\Bbbk)=q_{\mathrm{rat}}(\CC_{2,\mu}(\Bbbk))\)
for its image functor. The subscript distinguishes this map from
the pro-\(2\) map \(q\) below.
\end{definition}

On fundamental groups, filling \(0\) kills the loop
\(\GA\). Let \(\Pi\) be the resulting quotient and \(\pi\)
the quotient homomorphism:
\begin{equation}
 \Pi=\Gamma/\overline{\langle\!\langle\GA\rangle\!\rangle}
       =\widehat F^{(2)}(\GB,\GD),
 \qquad \pi:\Gamma\longrightarrow\Pi.
 \label{eq:quotient}
\end{equation}
Thus \(\Pi\) is the pro-\(2\) fundamental group of
\(\mathbb P^1\setminus\{1,-1,\infty\}\), with generating loops
\(\GB,\GD\). We also write \(e_B,e_D\) for their exponent
homomorphisms on \(\Pi\).

\begin{definition}[The integral forward maps]
Set
\begin{equation}
 q:\GTwo\longrightarrow\Pi,\qquad q(f,h)=\pi(h),
 \qquad \Ima=q(\GTwo),\qquad \cQ=q\circ\cE.
 \label{eq:integral-forward-maps}
\end{equation}
Thus, for \(f\in\GOne\),
\[
 \cQ(f)=\GB^{s_f}\pi(\Sch(f)),\qquad
 s_f=-e_B(\Sch(f)).
\]
For \(k\in\Pi\), the notation \(k(U,V)\) means
substitution \(\GB\mapsto U\), \(\GD\mapsto V\) in the
pro-\(2\) word \(k\).
\end{definition}

We call the inverse on the projected image the
\emph{Hirose--Sato map}. Its rational coefficient formula uses the
word operators of \cite{HiroseSato} and the mixed-pentagon
calculation of \cite{Hirose}.

\begin{proposition}[Ordinary recovery under three relations]
\label{prop:three-recovery}
Let \(\Psi(a,b,d)\) be group-like over a commutative
\(\mathbb Q\)-algebra, with \(c_b(\Psi)=0\). Assume
\begin{enumerate}[label=\textnormal{(\arabic*)},leftmargin=*]
\item Broadhurst duality \eqref{eq:broadhurst} for \(\Psi\),
\item the cyclotomic octagon \eqref{eq:octagon} for \(\Psi\),
\item ordinary duality and the hexagon
  \eqref{eq:ordinary-duality}--\eqref{eq:ordinary-hexagon}
  for the boundary value \(\Phi(u,v):=\Psi(u,v,0)\).
\end{enumerate}
Put \(g(u,v)=\Psi(0,u,v)\) and \(\rho=c_v(g)=c_d(\Psi)\).
Then
\begin{equation}
 \Phi(u,v)=e^{\rho v}g(v,u)g(u,v)^{-1}e^{-\rho u}.
 \label{eq:rational-ordinary-recovery}
\end{equation}
\end{proposition}

\begin{proof}
Put \(R=-u-v\). At \((a,b,d)=(0,u,v)\), the octagon reads
\[
 g(v,u)^{-1}e^{\mu v/2}\Psi(R,v,u)e^{\mu R/2}
       \Psi(R,u,v)^{-1}e^{\mu u/2}g(u,v)=1 .
\]
The reflected third argument in both middle factors is zero.
Thus (1) and the definition in (3) give
\[
 \Psi(R,v,u)=e^{-\rho v}\Phi(v,R)^{-1}e^{-\rho R},\qquad
 \Psi(R,u,v)^{-1}=e^{\rho R}\Phi(u,R)e^{\rho u}.
\]
The exponentials in \(\rho R\) cancel across \(e^{\mu R/2}\).
Applying the hexagon to \((v,u)\) and using ordinary duality gives
\[
 e^{\mu v/2}\Phi(v,R)^{-1}e^{\mu R/2}\Phi(u,R)e^{\mu u/2}
 =\Phi(u,v).
\]
The specialized octagon is now
\(g(v,u)^{-1}e^{-\rho v}\Phi(u,v)e^{\rho u}g(u,v)=1\),
which gives \eqref{eq:rational-ordinary-recovery}.
\end{proof}

The following result recovers the full cyclotomic component as well.

\begin{proposition}[The rational Hirose--Sato map]
\label{prop:rational-reconstruction}
For \((\Phi,\Psi)\in\CC_{2,\mu}\), let
\(g=q_{\mathrm{rat}}(\Phi,\Psi)\). Formula
\eqref{eq:rational-ordinary-recovery} recovers \(\Phi\).
The cyclotomic component \(\Psi\) is uniquely recovered by
\begin{equation}
 c_w(\Psi)=\lambda_\Phi\!\left(
  \operatorname{pr}_d\,g(M,-\partial_0)w\right),
 \qquad M=\partial_0+\partial_1,
 \label{eq:rational-reconstruction-main}
\end{equation}
for admissible words \(w\), followed by shuffle regularization
with \(c_a(\Psi)=c_b(\Psi)=0\). Here
\(\operatorname{pr}_d\) sets \(d=0\), and
\(\lambda_\Phi(w)=c_w(\Phi(a,b))\). The integral deletion
operators \(\partial_0,\partial_1\) and admissible words are
specified in Appendix~\ref{app:coefficient-proof}.
The morphism \(q_{\mathrm{rat}}\) is a closed immersion, and its
image functor is represented by its closed image.
These formulas define a natural inverse
\(\cS_{\mathrm{rat}}:\Ima_{\mathrm{rat}}\to\CC_{2,\mu}\):
\[
 \cS_{\mathrm{rat}}q_{\mathrm{rat}}=\id,
 \qquad q_{\mathrm{rat}}\cS_{\mathrm{rat}}=\id_{\Ima_{\mathrm{rat}}}.
\]
\end{proposition}

\begin{proof}
Equations \eqref{eq:ordinary-specialization} and
\eqref{eq:broadhurst} give the hypotheses of
Proposition~\ref{prop:three-recovery}. For the coefficients of
\(\Psi\), Appendix~\ref{app:rational-transport} applies the mixed
pentagon to obtain \eqref{eq:rational-reconstruction-main}.
The deletion operators
lower word length, so only finitely many terms contribute to a
fixed coefficient. Shuffle regularization recovers all remaining
coefficients. Each recovered coefficient is a polynomial over
\(\mathbb Q\) in finitely many coefficients of \(g\). The formulas
therefore recover every preimage naturally in the coefficient algebra.
The coefficients of \(\Phi\) and \(\Psi\) generate the coordinate
algebra of \(\CC_{2,\mu}\) over the chosen parameter base.
The polynomial recovery formulas show that the induced map on
coordinate algebras is surjective. Thus \(q_{\mathrm{rat}}\) is a
closed immersion and represents the stated image functor.
\end{proof}

We now give the Betti counterparts of the three relations.

\begin{definition}[The normalized covering graph and its relations]
\label{def:three-betti}
The normalized covering graph consists of the pairs \((f,h)\) with
\[
 f\in F'=\overline{[F,F]},\qquad
 s=-e_B(\Sch(f)),\qquad h=\GB^s\Sch(f).
\]
Thus \(j(h)=Y^sf\) and \(s=e_D(h)\).
Let \(\delta:\Gamma\to F\) fill the puncture \(-1\), sending
\((\GA,\GB,\GD)\) to \((X,Y,1)\). Define an involution
\(S\) of \(\Gamma\) by
\[
 S(\GA)=\GB,\qquad S(\GB)=\GA,\qquad
 S(\GD)=\GA^{-1}\GD^{-1}\GB^{-1}.
\]
Let \(\mathcal H_{123}\) be the subset of the normalized graph
satisfying the following three conditions:
\begin{enumerate}[label=\textnormal{(\arabic*)},leftmargin=*]
\item Betti Broadhurst duality:
\begin{equation}
 h=\GB^{-s}S(h)^{-1}\GA^{-s};
 \label{eq:betti-three-bd}
\end{equation}
\item the parameter-one cyclotomic octagon of
\cite[equation~(27)]{Enriquez};
\item the boundary identity \(\delta(h)=f\), together with
ordinary duality \eqref{eq:group-duality} and the three-cycle
\eqref{eq:group-threecycle} for \(f\).
\end{enumerate}
\end{definition}

We retain the octagon in the definition to match the rational
formulation of Proposition~\ref{prop:three-recovery}.
On this graph the two ordinary identities imply the octagon.
Write \(x=X,\ y=Y,\ z=(xy)^{-1}\),
\(v=f(x,y)\), and \(w=f(z,y)\).
In the inclusion \(j\), equation (27) has the four factors
\[
  (\Ad_x h_0)^{-1}(\Ad_{z^{-1}} h_1)\,h_1^{-1}h_0,
  \qquad h_0=y^s v,\quad h_1=y^s w .
\]
Indeed, \(z^2\) is the image of the fourth loop
\(\GB^{-1}\GA^{-1}\GD^{-1}\), and \(z^{-1}yz=xyx^{-1}\).
Cancelling the powers of \(y\) reduces the octagon to
\[
 (\Ad_{z^{-1}}w)w^{-1}=(\Ad_x v)v^{-1}.
\]
To see this, put \(c=f(z,x)\).
Ordinary duality and the three-cycle at \((x,y)\) give \(w=vc\).
At \((z,y)\) the third loop is \((zy)^{-1}=\Ad_z x=\Ad_{y^{-1}}x\).
The same identities give \(w=(\Ad_{y^{-1}}v)(\Ad_z c)\).
Substitute \(c=v^{-1}w\) and apply \(\Ad_{z^{-1}}\).

\begin{lemma}[The three Betti relations for an associator]
\label{lem:betti-three}
Every \((f,h)\in\GTwo\) belongs to \(\mathcal H_{123}\).
\end{lemma}

\begin{proof}
Choose a parameter-one reference pair \((\Phi_*,\Psi_*)\) over
\(\mathbb Q\) by \cite[Corollary~7.5]{Enriquez}, and use the
module action \(\Psi_*\mapsto e^{tb}\Psi_*\) to impose \(c_b(\Psi_*)=0\).
Work over \(\QQ\). Put \({\Psi_*^\sigma}=\Psi_*(a,d,b)\) and \(\rho=c_d(\Psi_*)\).
The action of \((f,h)\) gives another normalized pair
\((\Phi_*\beta_1(f),\Psi_*\beta(h))\), where
\begin{align*}
 \beta_1(X)&=e^a, &\beta_1(Y)&=\Phi_*^{-1}e^b\Phi_*,\\
 \beta(\GA)&=e^a, &\beta(\GB)&=\Psi_*^{-1}e^b\Psi_*,&
 \beta(\GD)&=e^{a/2}{\Psi_*^\sigma}^{-1}e^d\Psi_*^\sigma e^{-a/2}
\end{align*}
by \cite[equation~(28)]{Enriquez}. Both maps are faithful.
Indeed, the integral Magnus expansion of a free pro-\(2\) group
is faithful after reduction modulo two, and the indicated
substitutions for its Magnus variables have invertible linear
part over \(\QQ\).
By Proposition~\ref{prop:boundary-duality}, \(\Psi_*(a,b,0)=\Phi_*(a,b)\).
It follows on the generators that \(\beta|_{d=0}=\beta_1\delta\).
Applying the same boundary identity to the new pair gives
\(\beta_1(f)=\beta_1(\delta h)\), hence \(f=\delta h\).

For duality set \(e=-a-b-d\), \(\sigma_0(a,b,d)=(b,a,e)\),
\(\Psi_{*,\infty}=\Psi_*(e,b,d)\), \(\Psi_{*,\infty}^{\sigma}=\Psi_*(e,d,b)\), and \(L=e^{\rho b}\Psi_*\).
We claim \(\Ad_{L^{-1}}\sigma_0\beta=\beta S\).
On \(\GA,\GB\) this follows from
\(\sigma_0(\Psi_*)=e^{-\rho a}L^{-1}\). For \(\GD\), put
\(U=\Psi_*^{-1}e^{b/2}\Psi_{*,\infty}\) and \(V={\Psi_*^\sigma}^{-1}e^{d/2}\Psi_{*,\infty}^{\sigma}\).
The reference octagon and its \(b,d\) interchange read
\[
 V e^{e/2}U^{-1}=e^{-a/2}\Psi_*^{-1}e^{-b}\Psi_*,\qquad
 U e^{e/2}V^{-1}=e^{-a/2}{\Psi_*^\sigma}^{-1}e^{-d}{\Psi_*^\sigma} .
\]
Multiplying the second equality by the first gives
\(Ue^eU^{-1}=\beta(\GA^{-1}\GD^{-1}\GB^{-1})\).
Also, reference duality gives
\(L^{-1}e^{b/2}\sigma_0({\Psi_*^\sigma})^{-1}=Ue^{\rho e}\).
These two equalities prove the claim on \(\GD\).
Now \(c_d(\Psi_*\beta(h))=\rho+s\). Substituting \(\Psi_*\beta(h)\) into
\eqref{eq:broadhurst} and using the claim gives
\[
 \beta(h)=\Psi_*^{-1}e^{-s b}\Psi_*\,\beta(S(h))^{-1}e^{-s a}
         =\beta(\GB^{-s}S(h)^{-1}\GA^{-s}).
\]
Faithfulness gives the asserted Betti duality. The ordinary
relations hold by the definition of \(\GOne\), and the octagon
follows by the graph calculation above.
\end{proof}

\begin{theorem}[The integral Hirose--Sato map]
\label{thm:full-reconstruction}
On the projected set
\(\{\pi(h):(f,h)\in\mathcal H_{123}\}\), define
\begin{align}
 s(k)&=e_D(k),\notag\\
 r(k)&=Y^{s(k)}k(Y,X)k(X,Y)^{-1}X^{-s(k)},\label{eq:r-formula}\\
 h(k)&=\GB^{s(k)}\Sch(r(k)).\label{eq:s-formula}
\end{align}
For every \((f,h)\in\mathcal H_{123}\) with \(\pi(h)=k\),
one has \((r(k),h(k))=(f,h)\). In particular, on \(\Ima=q(\GTwo)\)
the map \(\cS(k)=(r(k),h(k))\) satisfies
\begin{equation}
 \cS q=\id_{\GTwo},\qquad q\cS=\id_{\Ima},\qquad
 r\cQ=\id_{\GOne},\qquad \cQ r=\id_{\Ima}.
 \label{eq:full-inverse-identities}
\end{equation}
\end{theorem}

\begin{proof}
Fill \(0\) in the upstairs octagon, that is, set
\((\GA,\GB,\GD)=(1,U,V)\), and put \(T=U^{-1}V^{-1}\).
In the free pro-\(2\) group on \(U,V\) the four factors give
\[
 k(V,U)^{-1}h(T,V,U)\,h(T,U,U^{-1}VU)^{-1}k(U,V)=1 .
\]
The third argument after reflection by \(S\) is \(1\) in both
middle factors. Using \eqref{eq:betti-three-bd} and \(\delta h=f\),
their product is \(V^{-s}f(V,T)^{-1}f(U,T)U^s\).
The three-cycle at \((V,U)\) and duality give
\(f(V,T)^{-1}f(U,T)=f(U,V)\). Solving for \(f\) gives
\eqref{eq:r-formula}. The graph identity \(j(h)=Y^sf\) then gives
\eqref{eq:s-formula}. Lemma~\ref{lem:betti-three} applies this
calculation to \(\GTwo\), and \(\cQ=q\cE\) gives the last two
identities.
\end{proof}

\begin{remark}[Schreier rewriting]
The evaluation of \eqref{eq:r-formula}--\eqref{eq:s-formula}
uses products, inverses, substitutions, \(\ZZ\)-powers, and
Schreier rewriting. In particular, it requires no division by \(2\).
For a finite word, Schreier rewriting is obtained by tracking the
parity of the \(X\)-exponent. A letter \(Y^{\pm1}\) contributes
\(\GB^{\pm1}\) at even parity and \(\GD^{\pm1}\) at odd parity.
A letter \(X\) contributes \(\GA\) exactly when the parity
changes from odd to even; \(X^{-1}\) contributes \(\GA^{-1}\)
exactly when it changes from even to odd. This extends continuously
to the pro-\(2\) groups.
\end{remark}

\begin{remark}[Comparison of the two reconstructions]
\label{sec:rational-comparison}
In the notation of Lemma~\ref{lem:betti-three}, put
\(G=\Psi_*(0,b,d)\). Setting \(a=0\) in \(\beta\) factors it through
\(\pi\), inducing a homomorphism \(\bar\beta\) on \(\Pi\).
For \(k\in\Ima\), the action of \(\cS(k)\) on the reference
pair has projection \(G\bar\beta(k)\). Uniqueness of rational
reconstruction thus gives
\begin{equation}
 \cS_{\mathrm{rat}}\bigl(G\bar\beta(k)\bigr)
 =\bigl(\Phi_*\beta_1(r(k)),\ \Psi_*\beta(h(k))\bigr).
 \label{eq:rational-integral-comparison}
\end{equation}
The reconstruction formulas are independent of the reference pair;
its choice is used to compare the two realizations.
\end{remark}

\section{Restricted tangent spaces and the Hirose--Sato map}
\label{sec:restricted-tangent}

We equip the deformation groups with the filtrations induced by
upstairs Zassenhaus order. The main result is
Theorem~\ref{thm:tangent-reconstruction}: cyclotomic projection is
strict and its restricted differential is injective, with inverse
given by an explicit coefficient operator.

The relevant group law is the deformation composition \(\star\).
Let
\(\tau=\Ad_X|_\Gamma\); thus
\[
 \tau(\GA)=\GA,\qquad \tau(\GB)=\GD,\qquad
 \tau(\GD)=\GA\GB\GA^{-1}.
\]
For \(h\in\Gamma\), the restriction \(\alpha_h\) of
\(\theta_{j(h)}\) to \(\Gamma\) is given by
\begin{align}
 \alpha_h(\GA)&=h\GA h^{-1},&
 \alpha_h(\GB)&=\GB,&
 \alpha_h(\GD)&=h\tau(h)^{-1}\GD\tau(h)h^{-1}.
 \label{eq:ihara-upstairs}
\end{align}
It is an automorphism, since it induces the identity on the
Frattini quotient of the free pro-\(2\) group. Put
\(h_1\star h_2=\alpha_{h_2}(h_1)h_2\).
The identity \(\alpha_{h_1\star h_2}=\alpha_{h_2}\alpha_{h_1}\)
follows from the substitutions on \(X,Y\). It gives associativity,
and \(\alpha_h^{-1}(h^{-1})\) is the inverse of \(h\) for \(\star\).
On \(\Pi\), let \(\sigma\) exchange \(\GB,\GD\), and define
\begin{equation}
 \bar\alpha_k(\GB)=\GB,\qquad
 \bar\alpha_k(\GD)=k\sigma(k)^{-1}\GD\sigma(k)k^{-1},
 \qquad k_1\star k_2=\bar\alpha_{k_2}(k_1)k_2.
 \label{eq:ihara-target}
\end{equation}
The identities \(\pi\tau=\sigma\pi\) and
\(\pi\alpha_h=\bar\alpha_{\pi h}\pi\) show that \(\pi\),
and therefore \(q\), are homomorphisms for these laws.
These formulas agree with the normalized composition of
Proposition~\ref{prop:enriquez}. Indeed, write
\(j(h_i)=Y^{s_i}f_i\). Since
\(\theta_{Y^{s_2}f_2}=\Ad_{Y^{s_2}}\circ\theta_{f_2}\),
\[
 j(h_1\star h_2)
 =Y^{s_1+s_2}\theta_{f_2}(f_1)f_2
 =Y^{s_1+s_2}(f_1\star f_2).
\]
The action \(\alpha_{h_2}\) is trivial on
\(\Gamma^{\mathrm{ab}}\), so \(e_B(h_1\star h_2)=0\)
and \(e_D(h_1\star h_2)=s_1+s_2\).
Uniqueness of the normalized lift identifies the composition laws.
Thus \(\Ima\) is a subgroup of \((\Pi,\star)\).

\begin{definition}[The induced filtrations]
Using the intrinsic Zassenhaus filtrations of
\eqref{eq:zassenhaus}, set
\begin{align}
 F_n\GTwo&=\{(f,h)\in\GTwo:h\in D_n\Gamma\},
 \label{eq:upstairs-filtration}\\
 F_n\Ima&=\Ima\cap D_n\Pi.
 \label{eq:image-filtration}
\end{align}
On the target \((\Pi,\star)\) use \(D_n\Pi\).
The generators \(\GA,\GB,\GD\) have degree one upstairs;
\(\GB,\GD\) have degree one on the quotient.
\end{definition}

\begin{lemma}
\label{lem:group-filtration}
These are separated restricted central filtrations for the
composition law \(\star\). Multiplication on their successive
quotients is addition of initial Magnus terms.
\end{lemma}

\begin{proof}
If \(h\in D_n\Gamma\), each conjugator in
\eqref{eq:ihara-upstairs} belongs to \(D_n\Gamma\). Consequently
\((\alpha_h-\id)(I_\Gamma^m)\subset I_\Gamma^{m+n}\).
The identical assertion holds for \(\bar\alpha_k\).
For \(U_h(v)=\alpha_h(v)h\), the operator \(U_h-\id\)
raises augmentation order by at least \(n\). The composition
identity for \(\alpha_h\) gives \(U_{h_1\star h_2}=U_{h_2}U_{h_1}\). Commutators of operators congruent to the
identity in orders \(m,n\) are congruent to the identity in order
\(m+n\), and \(U_h^2-\id=(U_h-\id)^2\) in characteristic two.
Evaluation at \(1\) gives the required commutator and square
bounds. In degree \(n\), \(h_1\star h_2\) has the sum of the two
initial coefficients. The same argument on \(\Pi\) proves the
target assertions. Separation follows from the Magnus embedding.
\end{proof}

\begin{definition}[Restricted tangent spaces]
The restricted tangent space of \(\GTwo\) is
\[
 \mathfrak g_2=\bigoplus_{n\ge1}F_n\GTwo/F_{n+1}\GTwo.
\]
The target is
\[
 \mathfrak p=\bigoplus_{n\ge1}D_n\Pi/D_{n+1}\Pi,
\]
using \(\star\) for its restricted Lie operations.
As a graded vector space, \(\mathfrak p\) identifies with
\(\Prim\FF\langle b,d\rangle\); upstairs the corresponding
ambient primitive space is \(\Prim\FF\langle a,b,d\rangle\).
Here \(a,b,d\) are the degree-one Magnus classes of
\(\GA,\GB,\GD\), and the tensor-algebra coproduct makes them
primitive.
\end{definition}

We first record how the covering changes the order of the ordinary
component. On the normalized covering graph it is
\(f=Y^{-e_D(h)}j(h)\).

\begin{lemma}[Vanishing of the ordinary initial component]
\label{lem:ordinary-order}
If \(h\in D_n\Gamma\) and \(e_B(h)=0\), then
\(Y^{-e_D(h)}j(h)\in D_{n+1}F\).
\end{lemma}

\begin{proof}
The covering inclusion sends the degree-one classes
\((a,b,d)\) to \((0,y,y)\). If \(\eta\) is the degree-\(n\)
coefficient of \(h-1\), then the degree-\(n\) coefficient of
\(j(h)-1\) is \(\eta(0,y,y)\). In the commutative polynomial
algebra the image of a primitive is a sum of pure \(2\)-powers.
Since \(e_B(h)=0\), we have
\(\eta(0,y,y)=\eta(0,0,y)\). The latter is also the
degree-\(n\) coefficient of \(Y^{e_D(h)}-1\), by the substitution
\((\GA,\GB,\GD)=(1,1,Y)\). Both group elements lie in \(D_nF\).
Their quotient therefore lies in \(D_{n+1}F\).
\end{proof}

For the mixed pentagon we use the Betti insertions
\(\delta_i\) corresponding to Enriquez's \(\partial_i\). They satisfy
\[
 \delta_3(h)\delta_1(h)=\delta_0(f)\delta_2(h)\delta_4(h).
\]
Identify the classes of the standard loops by
\(X_{0i}\mapsto a_i\), \(x_{ij}(0)\mapsto p_{ij}\),
\(x_{ij}(1)\mapsto m_{ij}\). The induced degree-one substitutions are
\begin{align*}
 \delta_3^{(1)}(a,b,d)
   &=(a_1,\ p_{12}+p_{13},\ m_{12}+m_{13}),\\
 \delta_1^{(1)}(a,b,d)
   &=(a_2+p_{12}+m_{12},\ p_{23},\ m_{23}),\\
 \delta_2^{(1)}(a,b,d)
   &=(a_1+a_2-p_{12}+m_{12},\ p_{13}+p_{23},\ m_{13}+m_{23}),\\
 \delta_4^{(1)}(a,b,d)
   &=(a_1,\ p_{12},\ m_{12}).
\end{align*}
Here the formulas are written over \(\mathbb Z\), before reduction
modulo two.
The ordinary insertion \(\delta_0\) sends \((X,Y)\) to
classes \((p_{12},p_{23})\). These follow directly from the
group formulas in \cite[equation~(24) and the proof of
Lemma~6.6]{Enriquez}: products give sums in degree one, and
conjugations disappear. For \(\delta_2\), the factor
\(x_{12}(0)^{-2}\) accounts for the difference from the
residue insertion; that difference vanishes over \(\FF\).

\begin{lemma}[Initial equations]
\label{sec:initial-equations}
For \((f,h)\in F_n\GTwo\), one has \(f\in D_{n+1}F\).
Let \(\phi(x,y)\) and \(\eta(a,b,d)\) be the degree-\(n\)
terms of \(f-1\) and \(h-1\), respectively. Thus \(\phi=0\).
The pair is primitive and satisfies
\begin{equation}
 \eta(a,0,0)=\eta(0,b,0)=0,\qquad
 \phi(x,0)=\phi(0,y)=0.
 \label{eq:axes}
\end{equation}
In the enveloping algebra of the holonomy Lie algebra from
Section~\ref{sec:preliminaries}, reduced modulo two, they satisfy
the initial mixed pentagon
\begin{align}
 &\eta(a_1,p_{12}+p_{13},m_{12}+m_{13})
   +\eta(a_2+p_{12}+m_{12},p_{23},m_{23})\notag\\
 &\quad=\phi(p_{12},p_{23})
   +\eta(a_1+a_2+p_{12}+m_{12},
            p_{13}+p_{23},m_{13}+m_{23})
   +\eta(a_1,p_{12},m_{12}).
 \label{eq:initial-pentagon}
\end{align}
The other initial equations needed below are
\begin{equation}
 \phi(a,b)=\eta(a,b,0).
 \label{eq:specialization}
\end{equation}
\begin{equation}
 \begin{aligned}
 \phi(x,y)&=\phi(y,x),\\
 \phi(x,y)+\phi(y,x+y)+\phi(x+y,x)&=0.
 \end{aligned}
 \label{eq:initial-threecycle}
\end{equation}
\begin{equation}
 \eta(a,b,d)+\eta(a,d,b)+\eta(a+b+d,d,b)
       +\eta(a+b+d,b,d)=0.
 \label{eq:initial-octagon}
\end{equation}

\end{lemma}

\begin{proof}
By Lemma~\ref{lem:ordinary-order}, \(f\in D_{n+1}F\).
The degree-\(n\) coefficients of group elements in \(D_n\)
are primitive. The exact zero exponent conditions on \(f,h\)
give the axis identities, including the vanishing of all pure
powers on those axes.

Modulo augmentation degree \(n+1\), products of elements
\(1+u\), with \(u\) of degree at least \(n\), become sums
of their degree-\(n\) terms. Substitutions act through their
degree-one parts. Apply this to the group mixed pentagon, using
the preceding substitutions reduced modulo two and
Lemma~\ref{lem:arrangement-graded-algebra}. This gives
\eqref{eq:initial-pentagon}.

The exact boundary identity \(\delta h=f\) of
Lemma~\ref{lem:betti-three} gives \eqref{eq:specialization} by
setting \(d=0\) in degree \(n\).

Taking initial terms in ordinary duality and the three-cycle
identity gives \eqref{eq:initial-threecycle}.
For the cyclotomic octagon, specialize
\cite[equation~(27)]{Enriquez} to \(N=2\) and scalar parameter
one. Its four occurrences of \(h\) have degree-one arguments
\[
 (a,d,b),\quad(-a-b-d,d,b),\quad
 (-a-b-d,b,d),\quad(a,b,d).
\]
The conjugations in these arguments do not change their initial
classes. Taking degree \(n\), and reducing signs modulo two,
gives \eqref{eq:initial-octagon}.
\end{proof}

We shall also consider the system
\eqref{eq:axes}--\eqref{eq:initial-octagon} independently of its
origin, allowing arbitrary homogeneous primitive pairs
\((\phi,\eta)\), without imposing \(\phi=0\) or assuming
that they arise from group solutions.

\begin{definition}[The linear Hirose--Sato operator]
Let \(\gamma(b,d)\) be homogeneous and primitive, with
\(\gamma(b,0)=0\), and put \(\ell_\gamma(t)=\gamma(0,t)\).
For degree \(n\), the polynomial \(\ell_\gamma\) is zero unless
\(n\) is a power of two. Define derivations \(P,E\) of
\(\FF\langle a,b,d\rangle\) by
\begin{equation}
 P(a)=P(b)=0,\quad P(d)=[a+b,d],\qquad
 E(a)=E(d)=[a,d],\quad E(b)=0,
 \label{eq:operators}
\end{equation}
and set \(Q=R_d+E\), where \(R_d(v)=vd\).
Put
\begin{align}
 \phi_\gamma(x,y)
   &=\gamma(y,x)+\gamma(x,y)+\ell_\gamma(x)+\ell_\gamma(y),
       \label{eq:linear-ordinary}\\
 \eta_\gamma(a,b,d)
   &=\phi_\gamma(a,b)+\gamma(P,Q)(1),
       \label{eq:linear-cyclotomic}\\
 \cR(\gamma)&=(\phi_\gamma,\eta_\gamma).\notag
\end{align}
Every term in these formulas has degree \(n\).

\end{definition}

\begin{proposition}[Coefficient reconstruction]
\label{prop:initial-reconstruction}
Every homogeneous primitive pair \((\phi,\eta)\) satisfying the
formal system \eqref{eq:axes}--\eqref{eq:initial-octagon} is recovered from
\(\gamma(b,d)=\eta(0,b,d)\) by
\eqref{eq:linear-ordinary}--\eqref{eq:linear-cyclotomic}.
\end{proposition}

\begin{proof}
The integral coefficient argument is given in
Appendix~\ref{app:coefficient-proof}. It uses the word evaluations
of the mixed pentagon and retains all primitive pure powers in the
regularization step.
\end{proof}

\begin{theorem}[Restricted differential of the Hirose--Sato map]
\label{thm:tangent-reconstruction}
The inverse \(\cS\) is filtered for
\eqref{eq:image-filtration} and \eqref{eq:upstairs-filtration}.
In fact,
\[
 q(F_n\GTwo)=F_n\Ima.
\]
The tangent projection is
\begin{equation}
 dq:\mathfrak g_2\longrightarrow\mathfrak p,
 \qquad (\phi,\eta)\longmapsto\eta(0,b,d),
 \label{eq:dq}
\end{equation}
and it is injective. Identify
\(\gr_F\Ima=\im(dq)\subset\mathfrak p\). Its inverse is
\begin{equation}
 d\cS:\im(dq)\longrightarrow\mathfrak g_2,
 \qquad d\cS=\cR|_{\im(dq)},\qquad
 d\cS\circ dq=\id_{\mathfrak g_2}.
 \label{eq:tangent-inverse}
\end{equation}
The image here consists of initial terms of actual solutions.
\end{theorem}

\begin{proof}
Take \(k\in\Ima\cap D_n\Pi\), and write
\(\cS(k)=(f,h)\). If \(h\) had a nonzero initial term
\(\eta_m\) in degree \(m<n\), the corresponding pair
\((\phi_m,\eta_m)\) would satisfy the initial equations above,
while \(\eta_m(0,b,d)=0\). Proposition~\ref{prop:initial-reconstruction}
would force \(\eta_m=0\), a contradiction. Thus
\(h\in D_n\Gamma\), and Lemma~\ref{lem:ordinary-order} gives
\(f\in D_{n+1}F\).

In degree \(n\), that proposition states that the initial pair of
\(\cS(k)\) is exactly \(\cR(\gamma_n)\), where \(\gamma_n\)
is the degree-\(n\) term of \(k-1\). Thus the reconstructed
class depends only on \(\gamma_n\). This proves that
\(d\cS\) is well defined and gives its formula.
Since \(q\) is a homomorphism for \(\star\), its inverse on
\(\Ima\) is also a homomorphism. The filtration statement
therefore applies to differences of representatives as well.
Taking initial terms in \(\cS q=\id\) proves
\eqref{eq:tangent-inverse}.
\end{proof}

\section[Comparison of the coefficient maps]{\texorpdfstring{Comparison with \cite{DD}}{Comparison of the coefficient maps}}
\label{sec:dd-comparison}
\label{sec:weighted-enriquez}

Transport the filtration of Section~\ref{sec:restricted-tangent}
through Enriquez's isomorphism:
\begin{equation}
 F_n\GOne=\cE^{-1}(F_n\GTwo),\qquad
 \mathfrak g_1=\gr_F\GOne.
 \label{eq:transported-source}
\end{equation}
Then \(d\cE:\mathfrak g_1\to\mathfrak g_2\) is an isomorphism.
Theorem~\ref{thm:tangent-reconstruction} immediately implies that
\(d\cQ=dq\circ d\cE\) is injective, with inverse on its image
\begin{equation}
 dr=(d\cE)^{-1}\circ\cR|_{\im(d\cQ)}.
 \label{eq:main}
\end{equation}
For \(\gamma\in\im(d\cQ)\), its reconstructed class satisfies
\[
 d\cE\bigl(dr(\gamma)\bigr)=d\cS(\gamma)=(0,\eta_\gamma),
\]
where the ordinary component vanishes by Lemma~\ref{lem:ordinary-order}.

\begin{definition}[Covering coordinates and weights]
Put \(\mathbf C=\GD\GB^{-1}=[X,Y]\). The elements
\(\GA,\GB,\mathbf C\) are another free basis of \(\Gamma\).
In mod-two Magnus coordinates write
\[
 x=X-1,\quad y=Y-1,\qquad
 a_0=\GA-1,\quad b_0=\GB-1,\quad c_0=\mathbf C-1.
\]
Define the decreasing covering filtration \(W_n\Gamma\) by
weighted order at least \(n\) in \(h-1\), with
\(|a_0|=|c_0|=2\) and \(|b_0|=1\). Downstairs \(W_nF\)
is ordinary Magnus order at least \(n\). The exact covering
substitutions are
\begin{equation}
 j(a_0)=x^2,\qquad j(b_0)=y,\qquad
 j(c_0)=[x,y](1+x)^{-1}(1+y)^{-1}.
 \label{eq:weighted-magnus-cover}
\end{equation}
Thus the graded covering homomorphism is
\begin{equation}
 \iota:\FF\langle A,B,C\rangle\longrightarrow
              \FF\langle x,y\rangle,
 \qquad (A,B,C)\longmapsto(x^2,y,[x,y]),
 \qquad |A|=|C|=2,\quad |B|=1.
 \label{eq:weighted-cover}
\end{equation}
The capitals here denote weighted initial forms, not the
degree-one puncture classes used in Section~\ref{sec:restricted-tangent}.
For the normalized deformation group, use the covering coordinate
\begin{equation}
 u(f,h)=\GB^{-e_D(h)}h=\Sch(f).
 \label{eq:unnormalized-covering-coordinate}
\end{equation}
The \(W\)-filtration on the two deformation groups in the following
statement is
\begin{equation}
 W_n\GOne=\GOne\cap W_nF,\qquad
 W_n\GTwo=\{(f,h):u(f,h)\in W_n\Gamma\}.
 \label{eq:weighted-deformation-filtration}
\end{equation}
In particular, the second condition measures \(u(f,h)\), rather
than \(h\), in the covering filtration.
\end{definition}

\begin{proposition}[Enriquez rewriting]
\label{prop:weighted-enriquez}
In the coordinates \eqref{eq:unnormalized-covering-coordinate},
the associated graded of Enriquez's isomorphism is
\begin{equation}
 \gr_W\cE:\psi(x,y)\longmapsto(\psi,\widetilde{\psi}),
 \qquad \widetilde{\psi}(x^2,y,[x,y])=\psi(x,y).
 \label{eq:weighted-enriquez}
\end{equation}
For an ordinary homogeneous Lie polynomial \(\psi\) of degree
\(n>1\), the unique preimage \(\widetilde{\psi}\) is an ordinary Lie
polynomial in \(A,B,C\), of weighted degree \(n\). It is computed
by writing \(\psi\) as a Lie expression in
\(e_j=\ad_x^j(y)\) and making the substitutions
\begin{equation}
 e_{2i}\longmapsto\ad_A^i(B),\qquad
 e_{2i+1}\longmapsto\ad_A^i(C),\qquad i\ge0.
 \label{eq:explicit-weighted-enriquez}
\end{equation}
\end{proposition}

\begin{proof}
The map \(\iota\) is injective: for lexicographic order \(x>y\),
the leading words \(xx,y,xy\) form a prefix code. Distinct words
in \(A,B,C\) consequently have distinct leading words after
substitution. Equation~\eqref{eq:weighted-magnus-cover} then gives
\(W_n\Gamma=j^{-1}(W_nF)\).

In the covering coordinate, \(\cE\) is exactly
\(f\mapsto(f,\Sch(f))\), so strictness of the covering embedding
gives \eqref{eq:weighted-enriquez}. These are group filtrations
for \(\star\): if \(f\in W_n\GOne\), then
\[
 (\theta_f-\id)(I_F^m)\subset I_F^{m+n}.
\]
As in Lemma~\ref{lem:group-filtration}, the composition formula
makes \(W_\bullet\GOne\) a separated restricted central
filtration. The equality \(W_n\Gamma=j^{-1}(W_nF)\) implies
\(W_n\GTwo=\cE(W_n\GOne)\), so the same holds upstairs.

The ideal generated by \(y\)
in the ordinary free Lie algebra is freely generated by the
\(e_j\), and contains every homogeneous component of degree
greater than one. In characteristic two,
\(\ad_x^2=\ad_{x^2}\), which verifies the substitutions
\eqref{eq:explicit-weighted-enriquez}. Injectivity of \(\iota\)
gives uniqueness.
\end{proof}

Let \(L_n\) be the reduction of the rational solution lattice of
Section~\ref{sec:preliminaries}, and assume \(n>2\).
For \(\bar\psi\in L_n\), write
\[
 \bar\psi(x,y)=\widetilde{\psi}(x^2,y,[x,y]).
\]
This is the substitution of \cite[Lemma~3.1]{DD}.
Our slot order is \((A,B,C)\), whereas \cite{DD} uses \((A,C,B)\).

\begin{definition}[The forward coefficient map]
Define the full forward coefficient map
\begin{equation}
 T_n:L_n\longrightarrow\Lie_{\FF}(B,C)_n,
 \qquad T_n(\bar\psi)=\widetilde{\psi}(0,B,C),
 \qquad |B|=1,\quad |C|=2.
 \label{eq:full-forward-coefficient-map}
\end{equation}
\end{definition}

\begin{proposition}[Forward coordinate comparison]
\label{prop:forward-coordinate-comparison}
Let \(\psi\) be an ordinary homogeneous Lie polynomial of
weight \(n>1\), and let \(\widetilde{\psi}\) be its unique covering
rewrite, so that \(\widetilde{\psi}(x^2,y,[x,y])=\psi(x,y)\).
In covering coefficient coordinates the forward composite is
\[
 \psi\longmapsto\widetilde{\psi}(A,B,C)\longmapsto\widetilde{\psi}(0,B,C).
\]
The second arrow is induced by \(\pi\). For an actual
\(f\in\GOne\), put \(s_f=-e_B(\Sch(f))\). The exact group
identity underlying the comparison is
\begin{equation}
 \GB^{-s_f}\cQ(f)=\pi(\Sch(f)).
 \label{eq:forward-normalization-comparison}
\end{equation}
This change of coordinates is reversible on the two images:
if \(k_0=\pi(\Sch(f))\), then
\[
 s_f=e_D(k_0),\qquad
 \cQ(f)=\GB^{e_D(k_0)}k_0.
\]
Equip \(\Pi\) with the weighted Magnus filtration in the free
basis \((\GB,\mathbf C)\), where \(\mathbf C=\GD\GB^{-1}\),
with weights \(1,2\), respectively. If an actual element
\(f\in W_n\GOne\) has weight-\(n\) initial term \(\psi\),
then the weight-\(n\) term of
\(\GB^{-s_f}\cQ(f)-1\) is \(\widetilde{\psi}(0,B,C)\), possibly zero.
Thus \(T_n\) is the restriction to \(L_n\) of the forward
coefficient formula in the unnormalized covering coordinate.
This restriction is defined for every element of \(L_n\),
independently of whether it lifts to a pro-\(2\) group solution.
\end{proposition}

\begin{proof}
In the free bases \((\GA,\GB,\mathbf C)\) upstairs and
\((\GB,\mathbf C)\) on the quotient, the map \(\pi\)
sends \(\GA\) to \(1\) and fixes \(\GB,\mathbf C\).
Its Magnus substitution is therefore exactly
\((a_0,b_0,c_0)\mapsto(0,b_0,c_0)\); on coefficient
polynomials it is \(\widetilde{\psi}\mapsto\widetilde{\psi}(0,B,C)\).
The equality \(h_f=\GB^{s_f}\Sch(f)\) gives
\eqref{eq:forward-normalization-comparison}. Since the
abelianization of \(f\) is zero, its lift \(u=\Sch(f)\)
satisfies \(e_B(u)+e_D(u)=0\). Thus
\(e_D(k_0)=e_D(u)=-e_B(u)=s_f\), which proves the
inverse coordinate formula.
The weighted assertion follows from
Proposition~\ref{prop:weighted-enriquez} and the displayed
Magnus substitution.
\end{proof}

\begin{definition}[The leading projection]
Let \(\mathcal B^rL_n\) consist of those elements whose covering
rewrite has at least \(r\) letters \(B\) in every monomial, and put
\[
 \mathcal L_{n,r}=\gr_{\mathcal B}^r L_n.
\]
A class is represented by its component \(\chi(A,B,C)\) with
exactly \(r\) letters \(B\). Its leading projection is
\begin{equation}
 q_{n,r}:\mathcal L_{n,r}\longrightarrow\Lie_{\FF}(B,C),
 \qquad f(B,C)=q_{n,r}(\chi)=\chi(0,B,C).
 \label{eq:dd-projection}
\end{equation}
Equivalently, if \(\operatorname{pr}_{B=r}\) denotes extraction
of the component with exactly \(r\) letters \(B\), then
\begin{equation}
 q_{n,r}([\bar\psi])
 =\operatorname{pr}_{B=r}T_n(\bar\psi),
 \qquad \bar\psi\in\mathcal B^rL_n.
 \label{eq:forward-leading-comparison}
\end{equation}
Indeed, setting \(A=0\) preserves the number of \(B\)'s,
so it commutes with this component extraction. This is the
additional leading operation in the map of \cite{DD}.
\end{definition}

\begin{remark}[Weights and symbolic length]
If \(\chi\ne0\), there is an integer \(m>0\) such that
\begin{equation}
 n=2m+r,\qquad \ell=m+r=\frac{n+r}{2}.
 \label{eq:dd-reindexing}
\end{equation}
Here \(m\) counts letters from \(\{A,C\}\), and \(\ell\) is the
upstairs word length of this component, counting each of \(A,B,C\)
once. The case \(m=0\) would be an ordinary
Lie polynomial in \(B\) alone, which vanishes for \(n>2\).

In the intrinsic upstairs Magnus grading,
\(\mathbf C=\GD\GB^{-1}\) has initial form \(b+d\).
Thus the change between its degree-one coefficient letters and
the covering-letter presentation is
\[
 (a,b,d)=(A,B,B+C).
\]
This linear change preserves symbolic length. At fixed covering
weight \(n\), \eqref{eq:dd-reindexing} identifies the
\(B\)-count with twice that length minus \(n\). This is a
comparison of homogeneous coefficient spaces; the group filtration
\(F\) and the covering filtration \(W\) are distinct filtrations.
\end{remark}

\begin{lemma}[An image bound from ordinary duality]
\label{lem:image-duality}
Let \(\psi(x,y)=\psi(y,x)\) be any homogeneous ordinary Lie
polynomial over \(\FF\) of weight \(n>2\). Suppose its covering
rewrite has \(B\)-count at least \(r\ge1\). Write
\(f=\operatorname{pr}_{B=r}\widetilde\psi(0,B,C)\).
Then \(f\) lies in the weight-\(n\), \(B\)-count-\(r\) part of
\begin{equation}
 J_+:=\Lie_{\FF}(b_k:k\ge1),
 \qquad b_k=\ad_C^k(B).
 \label{eq:dd-positive-image}
\end{equation}
Every element of \(J_+\) is invariant under \(B\mapsto B+C\).
\end{lemma}

\begin{proof}
In the elimination letters \(e_j=\ad_x^j y\), each \(e_{2i}\)
has one \(B\) and each \(e_{2i+1}\) has none. Their associative
words are independent, so each term in this expansion has at
least \(r\) even indices. We may assume \(n=2m+r\) with
\(m\ge1\), since otherwise \(f=0\). A word surviving in \(f\)
has downstairs bidegree \((m,m+r)\). Apply the even-index bound
to the component of bidegree \((m+r,m)\), which has \(m\)
elimination letters in each term. By symmetry, exchanging \(x,y\)
gives the component projecting to \(f\), with the same bound in
the exchanged alphabet \cite[Lemma~5.1]{DD}. Lemma~5.2 there
deduces \(f(B+s(C),C)=f(B,C)\) for every
polynomial \(s\), also after adjoining scalar indeterminates.
The Vandermonde argument of \cite[Proposition~5.3]{DD} then
excludes \(b_0=B\) from the free elimination alphabet
\((b_0,b_1,\ldots)\). These arguments use only ordinary Lie
membership and symmetry. Each \(b_k\) for \(k\ge1\) is fixed
under \(B\mapsto B+C\).
\end{proof}

\begin{remark}[The upper bound on the image and the bound on the source]
In particular \(\dim\im q_{n,r}\le\dim (J_+)_{n,r}\) follows
from ordinary duality and Lie membership alone. At \(r=0\) the image is zero in
weight \(n>2\), since it lies in the ordinary Lie algebra on \(C\).
To transfer either assertion to the source, the additional input is
\begin{equation}
 0\ne\chi\in\mathcal L_{n,r}\quad\Longrightarrow\quad
 \chi(0,B,C)\ne0 .
 \label{eq:leading-survival}
\end{equation}
Equivalently, some word of the leading \(B\)-component with
nonzero coefficient contains no \(A\). At fixed weight \(n\),
this says exactly that the shortest upstairs length survives
projection, by \eqref{eq:dd-reindexing}. 
The pentagon enters the source estimate of
\cite[Corollary~5.4]{DD} precisely in the proof of
\eqref{eq:leading-survival}. If \(\chi(0,B,C)=0\), the four-term
consequence of the pentagon gives \(\chi=P(A,B)\)
\cite[Lemma~4.3]{DD}. The reduced special identity and ordinary
symmetry then give \(P=0\) \cite[Proposition~4.4]{DD}.
Theorem~\ref{thm:dd-comparison} gives an explicit left inverse,
also derived in Appendix~\ref{app:dd-reconstruction}.
For the larger solution set \(\mathcal H_{123}\), leading
survival remains a separate question; see
Appendix~\ref{app:three-relations}.
\end{remark}

\begin{remark}[Deligne's rational motivic analogue]
For an element \(\sigma\) of the unipotent motivic Galois Lie algebra
of \(\mathrm{MT}(\mathbb Z[1/2])\), let \(F(e_0,e_1,e_{-1})\)
denote its infinitesimal action evaluated at the canonical path
\(1_{1,0}\) on the level-two curve.
Write \(F_\ell\) for the sum of its words of length \(\ell\).
The action respects the Tate grading and each loop letter has degree
one, so \(F_\ell\) is the image of the degree-\(\ell\) part of \(\sigma\).
Deligne's projection to \(L''\) sets \(e_1=0\) and is injective
on this image \cite[\S\S~3.2--3.4, 4.1 and Corollary~4.11(i)]{Deligne2010}.
Thus \(F_\ell\ne0\) implies \(F_\ell(e_0,0,e_{-1})\ne0\).
The coordinate \(u=(1-z)/(1+z)\) sends \(z=0\) to \(u=1\), so
\(e_1\) corresponds to our filled puncture. A word avoiding this
letter therefore occurs with nonzero coefficient already in the
shortest nonzero component, as in \eqref{eq:leading-survival}.
\end{remark}

\begin{definition}[The reconstruction operators in covering coordinates]
For an ordinary Lie polynomial \(f(B,C)\) with no linear part,
put
\begin{equation}
 \gamma(b,d)=f(b,b+d),\qquad
 \Slin(f)=\eta_\gamma(A,B,B+C),
 \label{eq:dd-linearized-input}
\end{equation}
where \(\eta_\gamma\) is the coefficient formula
\eqref{eq:linear-cyclotomic}, extended linearly across homogeneous
degrees. The condition \(\gamma(b,0)=f(b,b)=0\) permits this
substitution.

On \(\FF\langle A,B,C\rangle\)
define derivations \(\mathsf E,\mathsf V\) by
\begin{equation}
 \begin{gathered}
 \mathsf E(A)=\mathsf E(B)=0,\qquad
 \mathsf E(C)=[B,A+C],\\
 \mathsf V(A)=[A,C],\qquad
 \mathsf V(B)=\mathsf V(C)=0,
 \qquad \mathsf N=R_{A+C}+\mathsf V.
 \end{gathered}
 \label{eq:dd-reconstruction-operators}
\end{equation}
Set
\begin{equation}
 \rLead(f)=f(\mathsf E,\mathsf N)(1).
 \label{eq:dd-reconstruction}
\end{equation}
The first slot corresponds to \(B\), the second to \(C\).
\end{definition}

\begin{lemma}[Triangular form of the Hirose--Sato operator]
\label{lem:triangular-reconstruction}
Let \(\kappa\) be the derivation with
\(\kappa(A)=\kappa(B)=0\) and \(\kappa(C)=[A,C]\).
For every \(f\in J_+\),
\begin{equation}
 \Slin(f)=f(\mathsf E+\kappa,\mathsf N)(1).
 \label{eq:dd-triangular-formula}
\end{equation}
\end{lemma}

\begin{proof}
Define derivations
\(\partial,\Delta\) by
\begin{equation}
 \begin{gathered}
 \partial(A)=\partial(B)=0,\qquad
 \partial(C)=[A,B+C]+[B,C],\\
 \Delta(A)=\Delta(C)=0,\qquad \Delta(B)=[B,C],
 \qquad \mathsf M=L_{B+C}+\Delta.
 \end{gathered}
 \label{eq:dd-transport-operators}
\end{equation}
Under \((a,b,d)=(A,B,B+C)\), the operator \(P\) in
\eqref{eq:operators} becomes \(\partial\), and \(P+Q\)
becomes \(\mathsf M\). For example, the transported operator \(P\) kills \(A,B\)
and sends \(C\) to \([A+B,B+C]=[A,B+C]+[B,C]\).
The operator \(P+Q\) has multiplication part \(L_{B+C}\)
and derivation part \(\Delta\).
Substituting \(\gamma(b,d)=f(b,b+d)\) into
\eqref{eq:linear-cyclotomic} therefore gives the exact identity
\begin{equation}
 \Slin(f)=f(B,A+B)+f(A,A+B)+f(\partial,\mathsf M)(1).
 \label{eq:dd-full-linearized-formula}
\end{equation}
Here \(\ell_\gamma(t)=f(0,t)=0\). If \(f\in J_+\), its
invariance under \(B\mapsto B+C\) gives
\(f(B,A+B)=f(A,A+B)\), so the first two terms cancel.

The definitions give the exact operator identities
\begin{equation}
 \partial=\mathsf E+\kappa,\qquad
 \mathsf M=\mathsf N+R_A+L_B.
 \label{eq:dd-operator-decomposition}
\end{equation}
Both \(R_A\) and \(L_B\) commute with \(\partial\) and
\(\mathsf N\). Indeed, \(\partial(A)=\partial(B)=0\),
\(\mathsf V(B)=0\), and
\[
 [R_A,\mathsf N]=R_{[A,A+C]+\mathsf V(A)}=0.
\]
Thus \(R_A+L_B\) commutes with both entries of the Lie
expression. Every Lie monomial in \(f\) has length at least two,
so adding this commuting operator to the second entry leaves its
value unchanged.
This proves the stated identity.
\end{proof}

\begin{lemma}[The leading part of Willwacher's reconstruction]
\label{lem:willwacher-diagonal}
Let \(E_{\rm W},T\) be the derivations in \cite[equation~(13)]{Willwacher}.
On the ordered generators \((A,B,C)\) their values are
\[
 E_{\rm W}=([A,B],0,[A,B]),\qquad
 T=(0,[C,B],[A,B^2]).
\]
Put \(T_0=(0,[C,B],0)\), also a derivation. For every ordinary Lie
polynomial \(f(B,C)\) with zero linear part,
\begin{equation}
 \rLead(f)=f(E_{\rm W},L_C+T_0)(1)+f(E_{\rm W}+L_B,T_0)(1).
 \label{eq:willwacher-leading}
\end{equation}
\end{lemma}

\begin{proof}
Put \(M_0=L_C+T_0\). On the generators,
\(\mathsf E=E_{\rm W}+L_B+R_B\) and \(\mathsf N=M_0+R_A\).
The operator \(R_A\) commutes with \(\mathsf E,M_0\), since both
derivations \(\mathsf E,T_0\) kill \(A\). Removing \(R_A\)
therefore leaves the Lie expression unchanged. Next \(L_B\)
commutes with \(E_{\rm W}+R_B\) and with \(M_0\), the latter
because \(T_0(B)=[C,B]\). We may remove \(L_B\) from the first
slot. This gives \(f(E_{\rm W}+R_B,T_0+L_C)(1)\).

Derivations preserve each of the spaces of left and right
multiplications under commutators, and these two spaces commute.
Thus every expanded Lie term containing both types of
multiplication vanishes. Terms containing only derivations kill
\(1\). The two contributions that remain are
\(f(E_{\rm W},T_0+L_C)(1)\) and \(f(E_{\rm W}+R_B,T_0)(1)\).
In characteristic two \([R_u,R_v]=R_{[u,v]}\), and for any
derivation \(D\), \([D,R_u]=R_{D(u)}\). The same formulas hold
for left multiplication, and \(R_u(1)=L_u(1)\). This replaces
\(R_B\) by \(L_B\) in the second contribution and proves
\eqref{eq:willwacher-leading}.
\end{proof}

\begin{theorem}[Comparison of the reconstruction operators]
\label{thm:dd-comparison}
For every bihomogeneous \(f\in J_+\) of \(C\)-degree \(m\)
and \(B\)-degree \(r\), the polynomial \(\Slin(f)\) has symbolic
length \(\ell=m+r\), has no covering-weight component below
\(n=2m+r\), and satisfies
\begin{equation}
 \operatorname{pr}_{B=r}\Slin(f)
   =[\Slin(f)]_{W=n}
   =\rLead(f).
 \label{eq:dd-diagonal-comparison}
\end{equation}
Thus the covering-weight associated-graded coefficient operator
of \(\Slin\) on \(J_+\) is \(\rLead\).

For every \(\chi\in\mathcal L_{n,r}\),
\begin{equation}
 \rLead(q_{n,r}\chi)=\chi.
 \label{eq:dd-left-inverse}
\end{equation}
This gives an explicit inverse on the image of the injection
\(q_{n,r}\) of \cite[Proposition~4.4]{DD}.
\end{theorem}

\begin{proof}
By Lemma~\ref{lem:triangular-reconstruction}, it suffices to
compare the degrees of \(\mathsf E,\kappa,\mathsf N\).
Each raises symbolic length by one.
The operator \(\mathsf E\) raises \(B\)-degree by one,
whereas \(\kappa\) and \(\mathsf N\) preserve it.
There are exactly \(r\) first-slot occurrences in each word of
\(f\). To obtain output \(B\)-degree \(r\), one must choose
\(\mathsf E\) at every such occurrence. This gives
\(\rLead(f)\); every remaining term has smaller \(B\)-degree.
Since symbolic length is \(\ell\), covering weight is
\(2\ell-\deg_B\). The remaining terms therefore have weight
strictly greater than \(n\), proving
\eqref{eq:dd-diagonal-comparison}.

For \eqref{eq:dd-left-inverse}, choose \(\bar\psi\in\mathcal B^rL_n\)
representing \(\chi\) and let \(g=\widetilde\psi(0,B,C)\).
The reduced equations, including the special identity recorded
in Section~\ref{sec:preliminaries}, put \(\bar\psi\) in
Willwacher's \(\mathfrak{grt}_1(\FF)\). Thus
\cite[Proposition~3.2]{Willwacher}, obtained by evaluating the
infinitesimal pentagon, gives
\[
 \widetilde\psi=g(E_{\rm W},L_C+T)(1)+g(E_{\rm W}+L_B,T)(1).
\]
Each first-slot operator adds one \(B\) to every nonzero term.
Each second-slot operator preserves the count or adds two, the
latter precisely through \(T(C)=[A,B^2]\). The full output
has weight \(n\). To take its part of \(B\)-count \(r\), only
\(f=\operatorname{pr}_{B=r}g=q_{n,r}\chi\) contributes, and
one replaces \(T\) by \(T_0\). Lemma~\ref{lem:willwacher-diagonal}
identifies the result with \(\rLead(f)\), as claimed.
\end{proof}

\begin{corollary}[An even-weight coefficient vanishing]
For \(\bar\psi\in L_{2k}\) with \(k\ge2\),
\[
  c_{z_1x\,z_2x\cdots z_kx}(\bar\psi)=0
  \qquad (z_1,\ldots,z_k\in\{x,y\}).
\]
\end{corollary}

\begin{proof}
The component \(\chi_0=\widetilde\psi(A,0,C)\) of \(B\)-count zero
represents a class in \(\mathcal L_{2k,0}\). Its projection
\(\widetilde\psi(0,0,C)\) is an ordinary Lie polynomial in
\(C\) alone of degree greater than one, hence vanishes.
The injectivity in Theorem~\ref{thm:dd-comparison} gives \(\chi_0=0\).
Every word occurring in \(\widetilde\psi\) thus contains a \(B\).
Its \(B\)-count is even, since \(A,C\) have weight two.
Under \eqref{eq:weighted-cover}, the second \(B\) therefore places
a \(y\) in an even position, since each \(A\) or \(C\) contributes
two letters. This excludes every word in the statement.
\end{proof}
Ordinary duality alone allows the symmetric Lie polynomial
\(\ad_x^{2k-1}(y)+\ad_y^{2k-1}(x)\), whose coefficient of
\(yx^{2k-1}\) is \(1\) for \(k\ge2\). In odd weight, ordinary Lie membership and duality already imply the following coefficient identity.

\begin{proposition}[An odd-weight identity from ordinary duality]
Let \(m\ge1\) and let \(\psi\) be a homogeneous ordinary Lie polynomial
of weight \(2m+1\) over \(\FF\) satisfying \(\psi(x,y)=\psi(y,x)\).
Then
\[
  c_{(xy)^{m-1}xxy}(\psi)=c_{(xy)^m y}(\psi).
\]
In particular, \(c_{xyxxy}(\psi)=c_{xyxyy}(\psi)\) in weight five.
\end{proposition}

\begin{proof}
By symmetry it suffices to compare \(u=(xy)^{m-1}xxy\) and
\(v=(yx)^m x\). Write \(H=\widetilde\psi\) for the covering rewrite
of Proposition~\ref{prop:weighted-enriquez}. Under
\eqref{eq:weighted-cover}, a covering word with \(a\) letters \(A\)
and \(b\) letters \(B\) has \(2a-b\) more \(x\)'s than \(y\)'s.
Both \(u\) and \(v\) have excess one and contain exactly one pair
\(xx\). Thus only the component \(K\) of \(H\) with one \(A\),
one \(B\), and \(m-1\) copies of \(C\) contributes.

In \(u\) the pair \(xx\) fixes \(A\) and the final \(y\) fixes \(B\).
In \(v\) the final \(xx\) fixes \(A\), and \(B\) can occupy any
odd position of the alternating prefix. Each such choice determines
the other blocks uniquely. Hence
\[
  c_u(\psi)=c_{C^{m-1}AB}(K),\qquad
  c_v(\psi)=c_{C^m A}\bigl(K(A,C,C)\bigr).
\]
By Jacobi, Lie monomials with one \(B\) are spanned by iterated
adjoints of the other letters applied to \(B\). In the multidegree
of \(K\), these are
\[
  K_j=\ad_C^j[A,\ad_C^{m-1-j}B],\qquad 0\le j\le m-1.
\]
The unique word of \(K_j\) ending in \(B\) is \(C^jAC^{m-1-j}B\),
with coefficient one. On setting \(B=C\), only \(K_{m-1}\) survives.
It becomes \(\ad_C^m A\) over \(\FF\), with coefficient one on
\(C^m A\). Thus both functionals take the value \(\delta_{j,m-1}\)
on \(K_j\), proving the equality.
\end{proof}

\appendix
\section{The integral coefficient calculation}
\label{app:coefficient-proof}

\subsection{Word operators and the restricted reconstruction}

For a commutative ring \(R\), write
\(\mathcal W_R=R\langle a,b,d\rangle\), equipped with the
shuffle product, and define the admissible submodule by
\begin{equation}
 \mathcal W_R^{\mathrm{adm}}
 =R\oplus Rd\oplus
   \sum_{u\in\{a,d\},\ v\in\{b,d\}}u\mathcal W_Rv.
 \label{eq:admissible-words}
\end{equation}
Thus a nonempty admissible word is either \(d\), or begins in
\(a,d\) and ends in \(b,d\). Put
\(\mathcal W=\mathcal W_{\FF}\). Transposes below are taken
with respect to the word-coefficient pairing.

\begin{lemma}[Regularization without division]
\label{lem:regularization}
A coefficient functional associated to a primitive series over
\(\FF\) is determined by its values on
\(\mathcal W^{\mathrm{adm}}\) and on the pure powers
\(a^j,b^j\), \(j\ge1\).
\end{lemma}

\begin{proof}
The functional vanishes on every shuffle product of two nonempty
words. Write a word as \(b^i v a^j\), where either \(v\) is
empty or its first letter lies in \(\{a,d\}\) and its last
letter lies in \(\{b,d\}\).

If \(i>0\) and \(v\ne1\), the word \(b^i v a^j\)
appears with coefficient one in \(b^i\shuffle(va^j)\).
Every other word in this shuffle has a shorter initial
\(b\)-block. Induction removes the initial block. For a word
\(va^j\) with \(v\ne1\) and \(j>0\), the same argument
applied to \(v\shuffle a^j\) removes the terminal
\(a\)-block; every other term has fewer terminal \(a\)'s
and still begins in \(a,d\).

If \(v=1\) and \(i,j>0\), use
\(b^i\shuffle a^j\); all other terms have a shorter initial
\(b\)-block and fall under the preceding induction. The
remaining cases are admissible words, pure powers, and the empty
word, whose coefficient is zero for a primitive series. Each
induction uses a coefficient equal to one, so no division is needed.
\end{proof}

For \(w=u_1\cdots u_m\), set \(u_0=b\), \(u_{m+1}=a\),
\(e_0=a\), and \(e_1=b\). The following operators are defined
over \(\mathbb Z\), then reduced modulo two:
\begin{equation}
 \partial_c(w)=\sum_{j=1}^m
 \left(
 \delta_{\{u_{j-1},u_j\},\{d,e_c\}}
 -\delta_{\{u_j,u_{j+1}\},\{d,e_c\}}
 \right)u_1\cdots\widehat{u_j}\cdots u_m,
 \qquad c=0,1.
 \label{eq:deletion}
\end{equation}

The symbol \(\delta\) compares unordered pairs, and
\(\partial_c(1)=0\). Put \(M=\partial_0+\partial_1\).

\begin{lemma}[Deletion and duality]
\label{lem:deletion-duality}
The operators \(\partial_0,\partial_1\) lower word length by
one and preserve
\(\mathcal W_{\mathbb Z}^{\mathrm{adm}}\). Let
\(J_{\mathbb Z}\) be the word anti-automorphism defined by
\[
 J_{\mathbb Z}(a)=d-b,\qquad J_{\mathbb Z}(b)=d-a,\qquad
 J_{\mathbb Z}(d)=d.
\]
It preserves the admissible submodule and satisfies
\begin{equation}
 \partial_cJ_{\mathbb Z}=J_{\mathbb Z}\partial_c
 \quad\text{on }\mathcal W_{\mathbb Z}^{\mathrm{adm}}.
 \label{eq:word-duality-commutation}
\end{equation}
Moreover, if \(\operatorname{pr}_d\) sets \(d=0\), then
\begin{equation}
 \operatorname{pr}_dM=0,
 \qquad \eps M=0.
 \label{eq:epsilon-deletion}
\end{equation}
All these identities remain valid after reduction modulo two.
\end{lemma}

\begin{proof}
The duality identity is \cite[Theorem~10(3)]{HiroseSato},
after reversing words to use endpoints \((b,a)\). It holds
over \(\mathbb Z\), not merely over \(\mathbb Q\).
Preservation of admissibility follows directly from
\eqref{eq:deletion}: an endpoint deletion that would leave an
initial \(b\) or a terminal \(a\) has cancelling adjacent-edge
contributions. The anti-automorphism reverses words and interchanges
the allowed initial and terminal letters, so it too preserves
admissibility.

For \(\operatorname{pr}_dM\), only the deletion of the sole
\(d\) in a word can leave a word without \(d\). In that case
both neighbors, including the fixed endpoints, belong to
\(\{a,b\}\); their two contributions to \(M\) cancel.
Words without \(d\) are annihilated by both operators. This
proves \(\operatorname{pr}_dM=0\), and applying \(\eps\)
gives the second identity.
\end{proof}

\begin{lemma}[Word evaluation of the initial mixed pentagon]
\label{lem:word-evaluation}
Let \((\phi,\eta)\) be a homogeneous primitive pair satisfying \eqref{eq:axes},
\eqref{eq:initial-pentagon}, and \eqref{eq:specialization}, and
put \(k(a,d)=\eta(a,0,d)\). For every admissible word \(w\),
\begin{equation}
 c_w(\eta)=c_w(\phi(a,b))
       +\eps\bigl(k(\partial_0+\partial_1,\partial_1)w\bigr).
 \label{eq:initial-word-transport}
\end{equation}
\end{lemma}

\begin{proof}
We first construct the functional over \(\mathbb Z\). In the
arrangement coordinates used for \eqref{eq:initial-pentagon}, put
\[
 \xi_a=0,\qquad \xi_b=1,\qquad
 \xi_d=-\frac{z_3}{z_2},\qquad \zeta=\frac{z_1}{z_2}.
\]
For \(w=u_1\cdots u_m\), use the formal iterated-integral symbol
\[
 I_w=I(0;\xi_{u_m},\ldots,\xi_{u_1};\zeta).
\]
The reversal matches our coefficient convention. Pair
\(a_i,p_{ij},m_{ij}\) with the coefficients of the forms
\(d\log z_i\), \(d\log(z_i-z_j)\), and
\(d\log(z_i+z_j)\), respectively.
For a degree-one arrangement generator \(s\), define
\begin{equation}
 \begin{split}
 \delta_s I(t_0;t_1,\ldots,t_m;t_{m+1})
 ={}&\sum_{i=1}^m\sum_{\varepsilon=\pm1}
 \varepsilon\,\langle s,d\log(t_{i+\varepsilon}-t_i)\rangle\\
 &\hspace{1em}\cdot
 I(t_0;t_1,\ldots,\widehat t_i,\ldots,t_m;t_{m+1}),
 \end{split}
 \label{eq:integral-symbol-deletion}
\end{equation}
where \(d\log0=0\). Every nonconstant difference between
\(0,1,\xi_d,\zeta\) is, up to sign, one of
\[
 \frac{z_3}{z_2},\quad \frac{z_1}{z_2},\quad
 \frac{z_2+z_3}{z_2},\quad
 \frac{z_1-z_2}{z_2},\quad \frac{z_1+z_3}{z_2}.
\]
Thus all coefficients in \eqref{eq:integral-symbol-deletion} are
integers. Let \(\epsilon_I\) take value one on symbols with no
internal points and zero on the others. On the free tensor algebra
of arrangement generators, define
\begin{equation}
 \Lambda_w(s_1\cdots s_r)
   =\epsilon_I\bigl(\delta_{s_r}\cdots\delta_{s_1}I_w\bigr).
 \label{eq:integral-pairing}
\end{equation}
For \(r=0\), this means \(\Lambda_w(1)=\delta_{m,0}\).

The logarithmic forms above have integral coefficients in the
arrangement forms, so \(\Lambda_w\) is integer-valued on the
free tensor algebra over \(\mathbb Z\). After tensoring with
\(\mathbb Q\), it is the integrable pairing of
\cite[Definition~16 and Propositions~18--19]{Hirose}, with the
word order reversed as in \(I_w\). Hence it annihilates the
holonomy relation ideal over \(\mathbb Q\).

If \(v\) belongs to the integral relation ideal, then
\(\Lambda_w(v)\in\mathbb Z\) has zero image in \(\mathbb Q\),
so \(\Lambda_w(v)=0\). Thus the functional descends over
\(\mathbb Z\), and may then be reduced modulo two. This step
uses the integrality of the functional; it does not require a
characteristic-zero lift of \((\phi,\eta)\).

We now evaluate the five insertions. Write
\(L_a^t,L_b^t,L_d^t\) for deletion of the first letter, when it
is respectively \(a,b,d\). In coefficient-word order, the
three operators for \(\eta^{1,2,34}\) are
\(L_a^t,L_b^t,L_d^t\), while those for \(\eta^{1,2,3}\)
are \(L_a^t,L_b^t,0\). For \(\eta^{1,23,4}\), they are,
over \(\mathbb Z\),
\begin{equation}
 \delta_{a_1+a_2+p_{12}+m_{12}}=-M,\qquad
 \delta_{p_{13}+p_{23}}=0,\qquad
 \delta_{m_{13}+m_{23}}=\partial_1.
 \label{eq:middle-insertion-operators}
\end{equation}
For example, put \(s=m_{13}+m_{23}\). The only differences
in \eqref{eq:integral-symbol-deletion} with nonzero pairing
against \(s\) are, up to sign,
\[
 \xi_d-1=-\frac{z_2+z_3}{z_2},\qquad
 \xi_d-\zeta=-\frac{z_1+z_3}{z_2};
\]
both logarithmic differentials pair to one. Thus the endpoint
\(\zeta\) contributes as the letter \(b\). With
\(u_0=b\) and \(u_{m+1}=a\), reversing to coefficient-word
order gives
\[
 \delta_s(w)=\sum_{j=1}^m
 \bigl(\delta_{\{u_{j-1},u_j\},\{b,d\}}
       -\delta_{\{u_j,u_{j+1}\},\{b,d\}}\bigr)
 u_1\cdots\widehat{u_j}\cdots u_m
 =\partial_1(w).
\]
Likewise, \(a_1+a_2+p_{12}+m_{12}\) pairs to \(-1\) on
the three edges joining \(\xi_d\) to \(0,1,\zeta\), and
to zero on the other edges; this gives
\(-\partial_0-\partial_1=-M\). None of the logarithmic
differentials contains the forms dual to \(p_{13},p_{23}\),
giving the zero operator in \eqref{eq:middle-insertion-operators}.

For \(\eta^{12,3,4}\), the operators associated to
\(a_2+p_{12}+m_{12},p_{23},m_{23}\) preserve the span of
nonempty words ending in \(b\) or \(d\), and annihilate the
one-letter words \(b,d\). Thus every positive-length operator
product has zero empty-word evaluation on an admissible word.
For \(\phi^{2,3,4}\), \(\delta_{p_{12}}=L_b^t\) and
\(\delta_{p_{23}}=0\); its evaluation vanishes because
\(\phi(x,0)=0\).

Reduce these evaluations modulo two. In the order of
\eqref{eq:initial-pentagon}, they give
\begin{align*}
 \Lambda_w(\eta^{1,2,34})&=c_w(\eta),\\
 \Lambda_w(\eta^{12,3,4})&=0,\\
 \Lambda_w(\phi^{2,3,4})&=0,\\
 \Lambda_w(\eta^{1,23,4})
   &=\eps\bigl(k(M,\partial_1)w\bigr),\\
 \Lambda_w(\eta^{1,2,3})&=c_w(\phi(a,b)).
\end{align*}
The evaluation of \(\eta^{1,2,3}\) uses
\(\eta(a,b,0)=\phi(a,b)\). For \(\eta^{1,23,4}\), the
reversed operator order in
\eqref{eq:integral-pairing} gives \(\rev k\); since \(k\)
is primitive, \(\rev k=k\) in characteristic two. These are
the linear insertion evaluations underlying
\cite[Lemmas~22--23]{Hirose}. Applying \(\Lambda_w\) to
the initial mixed pentagon proves
\eqref{eq:initial-word-transport}.
\end{proof}

\begin{lemma}[Restricted duality]
\label{lem:restricted-duality}
Let \((\phi,\eta)\) satisfy the hypotheses of
Proposition~\ref{prop:initial-reconstruction}. Put
\(\gamma(b,d)=\eta(0,b,d)\) and
\(\ell(t)=\gamma(0,t)\). Then
\begin{equation}
 \eta(a,b,d)+\eta(b,a,a+b+d)=\ell(a)+\ell(b).
 \label{eq:restricted-duality}
\end{equation}
In particular, for \(k(a,d)=\eta(a,0,d)\),
\begin{equation}
 k(a,d)=\gamma(a,a+d)+\ell(a).
 \label{eq:k-gamma}
\end{equation}
\end{lemma}

\begin{proof}
Let \(J\) be the reduction of \(J_{\mathbb Z}\) modulo two.
Apply \eqref{eq:initial-word-transport} to an admissible word
\(w\) and to \(Jw\), extending the identity linearly.
Lemma~\ref{lem:deletion-duality} gives
\[
 \eps k(M,\partial_1)Jw
 =\eps Jk(M,\partial_1)w
 =\eps k(M,\partial_1)w.
\]
The boundary terms also agree. Indeed, setting \(d=0\) after
\(J\) exchanges \(a,b\) and reverses words; \(\phi\)
is symmetric in \(a,b\), and reversal fixes a primitive series
over \(\FF\). Thus \(c_w(\eta)=c_{Jw}(\eta)\).

Transposition of \(J\), followed by the same reversal identity,
gives
\[
 c_{Jw}(\eta)=c_w\bigl(\eta(b,a,a+b+d)\bigr).
\]
Consequently the left-hand side of
\eqref{eq:restricted-duality} vanishes on admissible words.
It is primitive, since linear substitutions preserve primitives.
The abelianization of a primitive series over \(\FF\) is a
sum of pure \(2\)-power monomials. The vanishing axes in
\eqref{eq:axes} therefore imply
\[
 \eta(t,0,t)=\eta(0,t,t)=\ell(t).
\]
It follows that the restrictions of the left-hand side to the
\(a\)- and \(b\)-axes are \(\ell(a)\) and \(\ell(b)\).
The right-hand side has exactly these restrictions and vanishes
on admissible words. Lemma~\ref{lem:regularization} proves
\eqref{eq:restricted-duality}. Setting \(b=0\) gives
\eqref{eq:k-gamma}.
\end{proof}

\begin{proof}[Proof of Proposition~\ref{prop:initial-reconstruction}]
We first recover the ordinary component. Set
\(\gamma(b,d)=\eta(0,b,d)\), \(\ell(t)=\gamma(0,t)\),
and \(z=x+y\). The specialization \(a=0,b=x,d=y\) of
\eqref{eq:initial-octagon} is
\[
 \gamma(x,y)+\gamma(y,x)+\eta(z,y,x)+\eta(z,x,y)=0.
\]
Restricted duality and \eqref{eq:specialization} give
\begin{align*}
 \eta(z,y,x)&=\phi(y,z)+\ell(z)+\ell(y),\\
 \eta(z,x,y)&=\phi(x,z)+\ell(z)+\ell(x).
\end{align*}
Since \(\phi(y,z)+\phi(x,z)=\phi(x,y)\) by ordinary
duality and the three-cycle identity, this yields
\[
 \phi(x,y)=\gamma(y,x)+\gamma(x,y)+\ell(x)+\ell(y)
 =\phi_\gamma(x,y).
\]

To recover the cyclotomic component, substitute
\eqref{eq:k-gamma} into \eqref{eq:initial-word-transport}.
Because \(M+\partial_1=\partial_0\) over \(\FF\) and
\(\eps\ell(M)=0\) by \eqref{eq:epsilon-deletion}, we obtain
\begin{equation}
 c_w(\eta)=c_w(\phi(a,b))+
 \eps\bigl(\gamma(M,\partial_0)w\bigr)
 \qquad(w\in\mathcal W^{\mathrm{adm}}).
 \label{eq:word-final}
\end{equation}
The transposes on the full word spaces are
\begin{equation}
 M^t=P,\qquad \partial_0^t=Q,
 \label{eq:deletion-transposes}
\end{equation}
with \(P,Q\) as in \eqref{eq:operators}. The interior terms
in \(\partial_0^t\) give the derivation
\(E(a)=E(d)=[a,d]\), \(E(b)=0\), and its endpoint term
gives \(R_d\); adding \(\partial_1^t\) leaves
\(P(a)=P(b)=0\), \(P(d)=[a+b,d]\).
Transposition reverses operator order, and \(\rev\gamma=\gamma\).
Thus \eqref{eq:word-final} says that \(\eta\) and
\[
 \eta_\gamma=\phi_\gamma(a,b)+\gamma(P,Q)(1)
\]
have the same admissible coefficients.

Both series are primitive. To check this for the second one,
\(P\) and \(Q\) are coderivations for the primitive-letter
tensor coproduct. In particular \(R_d\) is a coderivation,
since \(d\) is primitive. Coderivations are closed under
commutators and, in characteristic two, under squaring. Since
\(\gamma\) is a restricted Lie polynomial, \(\gamma(P,Q)\)
is a coderivation; its value at \(1\) is primitive.
Moreover \(\gamma\) has no constant term, and the deletion
operators annihilate every word without \(d\). Hence
\(\gamma(P,Q)(1)\) has zero coefficients on pure
\(a\)- and \(b\)-powers. The same is true of
\(\phi_\gamma(a,b)=\phi(a,b)\) by \eqref{eq:axes}.
The two series therefore agree on these powers as well.
Lemma~\ref{lem:regularization} gives
\(\eta=\eta_\gamma\), completing the reconstruction.
\end{proof}

\subsection{The rational coefficient formula}
\label{app:rational-transport}

Here \(\Bbbk\) is a commutative \(\mathbb Q\)-algebra.
The word spaces, projection \(\operatorname{pr}_d\), and
integral deletion operators are extended to \(\Bbbk\).

\begin{proof}[Proof of the coefficient formula
\eqref{eq:rational-reconstruction-main}]
Let \(g(b,d)=\Psi(0,b,d)\), as in
Proposition~\ref{prop:rational-reconstruction}. Put \(K(a,d)=\Psi(a,0,d)\). The same integral word
functional applied to the full mixed pentagon gives
\begin{equation}
 c_w(\Psi)=\lambda_\Phi\!\left(
       \operatorname{pr}_d\,
       (\rev K)(-M,\partial_1)w\right).
 \label{eq:rational-intermediate-transport}
\end{equation}
This is \cite[Lemmas~22--23]{Hirose} in the present word order.
Indeed, expand the two factors on the left in operator order.
The first removes a prefix of \(w\). The remaining suffix is
either empty or still ends in \(b,d\), so only the empty suffix
survives evaluation against the second factor. The result is
\(c_w(\Psi)\). The ordinary insertion acts trivially because
\(\Phi(x,0)=1\). On the other side, the middle insertion
acts by \((\rev K)(-M,\partial_1)\), by
\eqref{eq:middle-insertion-operators}, and the last insertion
evaluates the resulting word through
\(\lambda_\Phi\operatorname{pr}_d\).

Let \(\rho=c_d(g)\). Broadhurst duality
\eqref{eq:broadhurst}, specialized at its second variable equal
to zero, gives
\[
 K(a,d)=g(a,-a-d)^{-1}e^{-\rho a}.
\]
The group-like antipode identity is
\((\rev K)(-a,-d)=K(a,d)^{-1}\). Substituting operators
therefore gives
\begin{align*}
 (\rev K)(-M,\partial_1)
 &=K(M,-\partial_1)^{-1}\\
 &=e^{\rho M}g(M,-M+\partial_1)
 =e^{\rho M}g(M,-\partial_0).
\end{align*}
Since \(\operatorname{pr}_dM=0\),
\(\operatorname{pr}_d e^{\rho M}=\operatorname{pr}_d\).
Substitution into \eqref{eq:rational-intermediate-transport}
proves \eqref{eq:rational-reconstruction-main}. Each operator
lowers word length, so every displayed series acts by a finite sum
on a fixed word.

Finally, the coefficient functional of \(\Psi\) is a shuffle
character. Its pure \(a\)- and \(b\)-power coefficients vanish:
the one-variable specializations are exponentials with zero linear
term. The shuffle inductions used in
Lemma~\ref{lem:regularization}, now with the known product terms
retained, determine the remaining coefficients by induction on
length and on the boundary blocks.
\end{proof}

\section{An explicit reconstruction of the leading component}
\label{app:dd-reconstruction}

We derive \eqref{eq:dd-left-inverse} directly from the
deletion tests, special identity and image bound of \cite{DD}.
Retain \(n=2m+r>2\),
\(\ell=m+r\), \(\chi\in\mathcal L_{n,r}\), and
\(f=\chi(0,B,C)\) from Section~\ref{sec:dd-comparison}.
The operators \(\mathsf E,\mathsf V,\mathsf N\) and
\(\rLead\) are defined in
\eqref{eq:dd-reconstruction-operators}--\eqref{eq:dd-reconstruction}.
All Lie algebras in this appendix are ordinary Lie algebras over
\(\FF\).

\subsection{An exact identity from the four-term deletion tests}

We recall the four-term relation in our slot order. In the
arrangement notation of Section~\ref{sec:preliminaries}, retain
\(a_i\) and put
\[
 u_{ij}=p_{ij},\qquad v_{ij}=p_{ij}+m_{ij}.
\]
We use the integral holonomy presentation reduced modulo two.
Proposition~3.2 of \cite{DD} states that
\begin{align}
 &\chi(a_1,u_{12}+u_{13},v_{12}+v_{13})
   +\chi(a_2+v_{12},u_{23},v_{23})\notag\\
 &\quad+\chi(a_1+a_2+v_{12},u_{13}+u_{23},v_{13}+v_{23})
   +\chi(a_1,u_{12},v_{12})\equiv0
       \pmod{U_{\ell,<r}}.
 \label{dd:eq:four-term-relation}
\end{align}
Here \(U_{\ell,<r}\) is the span of products of \(\ell\) named
generators containing fewer than \(r\) symbols \(u\).

Define further derivations by
\begin{equation}
 \begin{gathered}
 \mathsf E_0(A)=\mathsf E_0(B)=0,\qquad
 \mathsf E_0(C)=[B,C],\\
 \mathsf V_0(A)=\mathsf V_0(C)=[A,C],\qquad
 \mathsf V_0(B)=0,\qquad \mathsf N_0=R_C+\mathsf V_0.
 \end{gathered}
 \label{dd:eq:untwisted-operators}
\end{equation}

\begin{lemma}[Reconstruction modulo the two-letter part]
\label{dd:lem:reconstruction-mod-AB}
For \(\chi\) and \(f\) as in \eqref{eq:dd-projection},
\begin{equation}
 \chi(A,B,C)+\chi(A,B,A)=\rLead(f).
 \label{dd:eq:reconstruction-mod-AB}
\end{equation}
\end{lemma}

\begin{proof}
Use the scalar deletion tests of \cite[Section~4.1]{DD}, whose
internal alphabet is \(\{e,0,p\}\), with fixed endpoints
\(0,w\). Identify the internal letters \(e,0,p\) with the
coefficient letters \(A,B,C\), respectively. Test
\eqref{dd:eq:four-term-relation} on all words of length \(\ell\)
containing exactly \(r\) zeros. Lemma~4.2 of \cite{DD} kills
the remainder \(U_{\ell,<r}\) and forces each \(B\)-slot to
delete a zero in every surviving contribution.

The effective transposes for the four summands, with the slots in
the order \((A,B,C)\), are
\begin{equation}
 \begin{aligned}
 \mathsf O_1&=(R_A,R_B,R_A),\\
 \mathsf O_2&=(R_A,\mathsf E_0,\mathsf N_0),\\
 \mathsf O_3&=(0,\mathsf E_0,\mathsf V_0),\\
 \mathsf O_4&=(R_A,R_B,R_{A+C}).
 \end{aligned}
 \label{dd:eq:four-transpose-rows}
\end{equation}
To verify the rows, let \(\lambda\) be the edge weights of a
slot in \cite[table~(4.6)]{DD}. The transpose inserts a letter
\(x\) between consecutive letters or endpoints \(u,v\) with
coefficient \(\lambda(u,x)+\lambda(x,v)\). In a \(B\)-slot
one retains only \(x=0\). For example, the last slot of row two
has weight one on \(ep,pw\). The edge \(ep\) contributes
\([A,C]\) on both \(A\) and \(C\), and the edge \(pw\) appends
\(C\). This gives \(\mathsf V_0+R_C\). In row three the added
\(v_{13}\) cancels the weight on \(pw\), leaving \(\mathsf V_0\).

Reversal fixes \(\chi\), so transposition does not change the
Lie polynomial. Since
\([R_A,\mathsf E_0]=[R_A,\mathsf N_0]=0\),
any Lie monomial of length at least two involving the first slot
vanishes after this substitution. Thus the second summand contributes
\(f(\mathsf E_0,\mathsf N_0)(1)\). The third summand contributes zero,
since its remaining operators are derivations and annihilate \(1\).
The remaining summands
read the indicated coefficients, giving
\begin{equation}
 \chi(A,B,A)+f(\mathsf E_0,\mathsf N_0)(1)+\chi(A,B,A+C)=0.
 \label{dd:eq:untwisted-reconstruction}
\end{equation}
Every term has symbolic length \(\ell\) and \(B\)-degree \(r\),
so the tested words determine this entire coefficient identity.

Let \(\theta\) fix \(A,B\) and send \(C\) to \(A+C\).
It is an involution, and evaluation on generators gives
\[
 \theta \mathsf E_0\theta=\mathsf E,\qquad \theta \mathsf N_0\theta=\mathsf N.
\]
Applying \(\theta\) to \eqref{dd:eq:untwisted-reconstruction}
proves \eqref{dd:eq:reconstruction-mod-AB}.
\end{proof}

Write
\begin{equation}
 P(A,B)=\chi(A,B,A),\qquad \chi=\rLead(f)+P.
 \label{dd:eq:remainder}
\end{equation}
The remaining polynomial \(P\) is an ordinary Lie polynomial
with \(m\) letters \(A\) and \(r\) letters \(B\).

\subsection{The leading special identity}

Define a derivation \(\mathsf D\) and an involutive automorphism
\(\sigma_{\mathrm{sp}}\) of \(\FF\langle A,B,C\rangle\) by
\begin{equation}
 \mathsf D(A)=\mathsf D(C)=0,\quad \mathsf D(B)=C,\qquad
 \sigma_{\mathrm{sp}}(A)=A+C,\quad\sigma_{\mathrm{sp}}(B)=B,\quad\sigma_{\mathrm{sp}}(C)=C.
 \label{dd:eq:special-operators}
\end{equation}
They commute.

\begin{lemma}
\label{dd:lem:special-source}
If \(r\geq1\) and \(\chi\in \mathcal L_{n,r}\), then
\begin{equation}
 \mathsf D(\chi+\sigma_{\mathrm{sp}}\chi)=0.
 \label{dd:eq:leading-special-identity}
\end{equation}
\end{lemma}

\begin{proof}
Use the reduced special identity of \cite[Proposition~2.5]{DD}:
\[
 [x,\bar\psi(x,y)+\bar\psi(x+y,y)]
 +[y,\bar\psi(x+y,y)]=0.
\]
In covering coordinates, \(x\mapsto x+y\) sends
\(A\mapsto A+C+B^2\) and fixes \(B,C\). The part of
\(\ad_x\) lowering \(B\)-count is \(\mathsf D\); its
other part and \(\ad_y\) raise that count. Taking count
\(r-1\) gives the assertion, as in
\cite[proof of Proposition~4.4]{DD}.
\end{proof}

\begin{lemma}
\label{dd:lem:special-reconstruction}
For every \(f\in J_+\), the element \(h=\rLead(f)\) is an
ordinary Lie polynomial and
\begin{equation}
 \mathsf D(h+\sigma_{\mathrm{sp}} h)=0.
 \label{dd:eq:special-reconstruction}
\end{equation}
If \(f\) has \(m\) letters \(C\) and \(r\) letters \(B\),
then \(h\) has symbolic length \(m+r\) and \(B\)-degree \(r\).
\end{lemma}

\begin{proof}
Operators \(R_a+\beta\), with \(a\) an ordinary Lie polynomial
and \(\beta\) a Lie derivation, are closed under commutators.
Consequently
\begin{equation}
 \rho(f):=f(\mathsf E,\mathsf N)=R_h+\beta_f,
 \qquad h=\rho(f)(1),\quad \beta_f=f(\mathsf E,\mathsf V).
 \label{dd:eq:operator-decomposition}
\end{equation}
This gives ordinary Lie membership. The length and \(B\)-degree
follow from the operator degrees in Theorem~\ref{thm:dd-comparison}.

Direct calculation gives
\[
 \sigma_{\mathrm{sp}} \mathsf V\sigma_{\mathrm{sp}}=\mathsf V,\qquad
 \sigma_{\mathrm{sp}} \mathsf E\sigma_{\mathrm{sp}}=\mathsf E+\ad_B.
\]
The derivation \(\ad_B\) commutes with \(\mathsf E,\mathsf V\), since both
annihilate \(B\). Every Lie monomial in \(f\in J_+\) has
symbolic length at least two. Adding this commuting derivation
to its first slot has no effect, so
\begin{equation}
 \sigma_{\mathrm{sp}}\beta_f\sigma_{\mathrm{sp}}=\beta_f.
 \label{dd:eq:invariant-derivation}
\end{equation}

Put \(\mathsf K=[\mathsf D,\mathsf E]\). Checking the generators gives
\begin{gather*}
 \mathsf K(A)=0,\quad \mathsf K(B)=[B,A+C],\quad \mathsf K(C)=[A,C],\\
 [\mathsf D,\mathsf N]=0,\qquad [\mathsf V,\mathsf K]=\ad_{[A,C]},\qquad
 [\mathsf N,\mathsf K]=L_{[A,C]}.
\end{gather*}
Moreover, \([\mathsf N,L_a]=L_{\mathsf V(a)}\), and \(\mathsf V=\ad_C\) on the
subalgebra in \(A,C\). It follows that
\begin{equation}
 [\mathsf D,\rho(b_k)]=\ad_{\mathsf N}^k\mathsf K=L_{\ad_C^k A}
 \qquad(k\geq1).
 \label{dd:eq:generator-commutator}
\end{equation}
The property that \([\mathsf D,\rho(f)]\) is left multiplication is
closed under Lie brackets in \(f\), by Jacobi and
\([L_a,R_h+\beta]=L_{\beta(a)}\) in characteristic two.
Equation~\eqref{dd:eq:generator-commutator} therefore implies it
for all \(f\in J_+\). Evaluation on \(1\) identifies the
multiplication operator, giving
\[
 [\mathsf D,\rho(f)]=L_{\mathsf D h},\qquad
 [\mathsf D,\beta_f]=L_{\mathsf D h}+R_{\mathsf D h}=\ad_{\mathsf D h}.
\]
Both \(\mathsf D\) and \(\beta_f\) commute with \(\sigma_{\mathrm{sp}}\).
Hence \(\ad_{\mathsf D h+\sigma_{\mathrm{sp}}(\mathsf D h)}=0\). The center of the free
associative algebra on \(A,B,C\) consists of scalars, whereas
\(\mathsf D h+\sigma_{\mathrm{sp}}(\mathsf D h)\) has zero constant term. It is therefore
zero. Since \(\mathsf D\sigma_{\mathrm{sp}}=\sigma_{\mathrm{sp}} \mathsf D\), this is
\eqref{dd:eq:special-reconstruction}.
\end{proof}

\subsection{Eliminating the remainder}

\begin{proof}[Proof of \eqref{eq:dd-left-inverse}]
By \eqref{dd:eq:remainder}, it remains to show \(P=0\).
For \(r=0\), both \(f\) and \(P\) are ordinary Lie
polynomials in a single letter of weight two, so they vanish
when \(n>2\). Suppose \(r\ge1\). The image bound gives
\(f\in J_+\), and the preceding two lemmas yield
\begin{equation}
 \mathsf D\bigl(P(A+C,B)+P(A,B)\bigr)=0.
 \label{dd:eq:remainder-special}
\end{equation}
In a target word with two \(C\)'s, one came from a substituted
\(A\) and the other from a differentiated \(B\). The two ways
to assign these positions show that exchanging any \(A\) and
\(B\) leaves the coefficient in \(P\) unchanged. All coefficients
of bidegree \((m,r)\) therefore equal a scalar \(\lambda\).
If \(m\ge2\), take a word with three \(C\)'s, \(m-2\) letters \(A\)
and \(r-1\) letters \(B\). There are three choices for the
differentiated \(B\), so its coefficient in
\eqref{dd:eq:remainder-special} is \(3\lambda=\lambda=0\)
\cite[equations~(4.9)--(4.10)]{DD}.

Suppose \(m=1\) and \(r>1\). Each generator of \(J_+\) has
at least as many \(C\)'s as \(B\)'s, so \(f=0\) and
\(\chi=P(A,B)\). Since \(n=r+2\), the only possible higher
\(B\)-count term of the covering rewrite is a pure \(B\)-power.
It vanishes by ordinary Lie membership. Thus
\(\bar\psi=P(x^2,y)\) has bidegree \((2,r)\).
If \(r>2\), symmetry interchanges two distinct bidegrees, so
\(\bar\psi=0\) and \(P=0\) by injectivity of \(\iota\).
If \(r=2\), the ordinary Lie component is spanned by
\([B,[B,A]]=BBA+ABB\). Its coefficient on \(BAB\) is zero,
whereas \eqref{dd:eq:remainder-special} forces the three
coefficients to be equal. Hence \(P=0\) also in this case.
For \(m=r=1\), write \(f=\lambda[B,C]\) and
\(P=\nu[A,B]\). The reconstruction modulo \(P\) reads
\[
 \chi=(\lambda+\nu)[A,B]+\lambda[B,C].
\]
At weight three there is no higher \(B\)-count component:
it would be an ordinary Lie polynomial in \(B\) alone.
Thus \(\chi(x^2,y,[x,y])\) is symmetric under exchanging
\(x,y\). This exchanges the two displayed brackets after
substitution, so \(\nu=0\).
This proves the reconstruction identity.
\end{proof}

\section{Lifting on the graph of the three Betti relations}
\label{app:three-relations}

We study the solution set \(\mathcal H_{123}\) defined in
Section~\ref{sec:reconstruction}, without imposing the mixed
pentagon. We first parametrize the normalized covering graph with
its boundary condition, then impose Broadhurst duality and the two
ordinary identities at successive finite precisions. This gives
an obstruction-theoretic description of the realized initial
coefficients, rather than an explicit determination in every degree.
The octagon is automatic on this graph once the ordinary identities
hold.

Let \(T_{123,n}\) be the set of degree-\(n\) coefficients of
\(h-1\) for \((f,h)\in\mathcal H_{123}\) with \(h\in D_n\Gamma\),
including zero. These are the coefficients realized by exact
solutions. The following description uses only the solution set
and the filtration of the ambient free pro-\(2\) group.

Let \(K=\ker\delta\), let \(i:F\to\Gamma\) send \(X,Y\) to
\(\GA,\GB\), and put \(\alpha=j\circ i:F\to F\). Thus \(\alpha(X)=X^2\)
and \(\alpha(Y)=Y\).
\begin{proposition}[Boundary and lifting for the three relations]
\label{prop:three-tangent}
The normalized graph with the boundary condition \(\delta(h)=f\)
is parametrized by \(k\in K\). Its parameter and inverse are
\begin{equation}
 \begin{gathered}
 k=h\,i(f)^{-1},\qquad s=e_D(k),\qquad w=Y^{-s}j(k),\\
 f_k=\prod_{j\ge0}\alpha^j(w),\qquad h_k=k\,i(f_k),
 \end{gathered}
 \label{eq:three-boundary-param}
\end{equation}
with factors ordered by increasing \(j\). Its degree-\(n\)
initial coefficients form exactly
\[
 J_n=\ker\bigl(\Prim_n\FF\langle a,b,d\rangle
       \xrightarrow{\,d=0\,}\Prim_n\FF\langle a,b\rangle\bigr).
\]
Here \(\Prim_n\) denotes the homogeneous degree-\(n\) primitive
subspace.
For \(M\ge n\), run through \(k\in K\cap D_n\Gamma\) modulo
\(K\cap D_{M+1}\Gamma\). Let \(E_{n,M}\subset J_n\) consist
of the degree-\(n\) coefficients for which \(h_k\) satisfies
\eqref{eq:betti-three-bd} modulo \(D_{M+1}\Gamma\) and \(f_k\)
satisfies \eqref{eq:group-duality} and \eqref{eq:group-threecycle}
modulo \(D_{M+1}F\). Then
\begin{equation}
 T_{123,n}=\bigcap_{M\ge n}E_{n,M}.
 \label{eq:three-realization}
\end{equation}
In passing from precision \(m-1\) to \(m\), a correction
\(\xi\in J_m\) changes the degree-\(m\) Broadhurst defect by
\begin{equation}
 \operatorname{BD}_m(\xi)=\xi(a,b,d)+\xi(b,a,a+b+d)
                  +\ell_\xi(a)+\ell_\xi(b),\qquad
 \ell_\xi(t)=\xi(0,0,t).
 \label{eq:three-bd-layer}
\end{equation}
It leaves the two ordinary defects fixed through degree \(m\).
\end{proposition}

\begin{proof}
On the boundary write \(h=k\,i(f)\). The graph identity becomes
\(f=w\alpha(f)\). The normal generator \(\GD\) of \(K\)
maps to a conjugate of \(Y\). Hence the abelianization of \(j(k)\)
is \(Y^s\), and \(w\in F'\). In Magnus coordinates \(\alpha^j\) substitutes
\(x^{2^j}\) for \(x=X-1\) and fixes \(y=Y-1\).
Every monomial of \(w-1\) contains \(x\), so
\(\alpha^j(w)-1\) has order at least \(2^j\), uniformly in
\(k\). Thus only finitely many factors contribute through any
fixed degree. The product in \eqref{eq:three-boundary-param}
converges and depends continuously on \(k\). Iterating the same
equation on \(F'\) proves uniqueness. By construction \(f_k\in F'\),
\(\delta(h_k)=f_k\) and \(j(h_k)=Y^sf_k\). The first two facts
give \(e_A(h_k)=e_B(h_k)=0\), and \(e_D(h_k)=e_D(k)=s\).

If \(k\in K\cap D_n\Gamma\), then \(e_B(k)=0\), so
Lemma~\ref{lem:ordinary-order} gives \(w\in D_{n+1}F\). Thus \(f_k\in D_{n+1}F\) and
\(h_k\equiv k\bmod D_{n+1}\Gamma\). Conversely, the lemma applied to
\(h\in D_n\Gamma\) gives \(f\in D_{n+1}F\), so
\(k=h\,i(f)^{-1}\equiv h\bmod D_{n+1}\Gamma\). Every \(\eta\in J_n\) is
realized on the boundary. Indeed, lift it to \(v\in D_n\Gamma\)
by the free restricted Lie description of the Zassenhaus graded. Then
\(\delta(v)\in D_{n+1}F\), so \(v\,i(\delta v)^{-1}\) belongs to
\(K\) and still has initial coefficient \(\eta\).

To compute the remaining conditions, replace \(k\) by \(k l\)
with \(l\in K\cap D_m\Gamma\) of degree-\(m\) coefficient
\(\xi\). In \(w\) the change from \(j(l)\) cancels in degree \(m\)
against the change of \(Y^{-s}\). Conjugation by an earlier word
acts trivially on that graded layer. Since \(\alpha\) preserves
the Zassenhaus filtration, each factor \(\alpha^j(w)\), and hence
the convergent product defining \(f_k\), is fixed modulo
\(D_{m+1}F\). The degree-\(m\) change of \(h_k\) is therefore
\(\xi\).
Write Broadhurst duality \eqref{eq:betti-three-bd} as \(pS(p)=1\), where
\(p=\GB^s h_k\). The change of \(p\) in this layer is
\(\xi+\ell_\xi(b)\). The linear substitution of \(S\) is
\((a,b,d)\mapsto(b,a,a+b+d)\), giving
\eqref{eq:three-bd-layer}. At each step one tests the ordinary
defects and solves this linear equation against the current
Broadhurst defect. All computations through degree \(M\) factor
through finite quotients. For a fixed \(\eta\in\bigcap_M E_{n,M}\),
the sets of lifts in \(K\cap D_n\Gamma\) passing each precision
are nested, closed and nonempty. Compactness gives an exact lift.
The octagon holds by the graph calculation in
Section~\ref{sec:reconstruction}. This proves
\eqref{eq:three-realization}.
\end{proof}

In the first layer one has
\[
 E_{n,n}=\ker\bigl(\operatorname{BD}_n:J_n\to
             \Prim_n\FF\langle a,b,d\rangle\bigr).
\]
Indeed, every element of \(J_n\) has a boundary lift and its
ordinary component vanishes through degree \(n\). Thus the
Broadhurst defect is the sole test in that layer. The intersection
in \eqref{eq:three-realization} selects the classes which also
pass all higher lifting conditions.

\begin{remark}[The first lifting obstruction]
The first-layer test need not suffice for realization. In degree
one, \(J_1=\FF d\) and \(\operatorname{BD}_1(d)=0\), so
\(E_{1,1}=J_1\). Take the boundary parameter \(k=\GD\).
Then \(s=1\) and
\[
 w=Y^{-1}XYX^{-1},\qquad
 f_k=\prod_{j\ge0}\alpha^j(w).
\]
The degree-two Magnus term of \(w-1\) is \([x,y]\), whereas
\(\alpha^j(w)-1\) has order at least three for \(j\ge1\).
Thus the degree-two term of \(f_k-1\) is \([x,y]\).
For \(\phi(x,y)=[x,y]\), the three-cycle defect is
\[
 \phi(x+y,x)+\phi(y,x+y)+\phi(x,y)=[x,y]\ne0.
\]
Every other boundary parameter with initial coefficient \(d\)
differs from \(\GD\) by an element of \(K\cap D_2\Gamma\).
Proposition~\ref{prop:three-tangent} shows that such a correction
leaves the ordinary component fixed through degree two.
Consequently \(d\notin E_{1,2}\), and
\[
 E_{1,1}=\FF d,\qquad E_{1,2}=T_{123,1}=\{0\}.
\]
\end{remark}

\begin{remark}[The reflected boundary on the covering graph]
For \((f,h)\in\mathcal H_{123}\) with \(f\ne1\), put
\(s=e_D(h)\) and let \(\widetilde\psi\) be the
first covering-weight component of \(u=\Sch(f)\), of weight \(n\).
Write \(\chi\) for its component of least \(B\)-count \(r\). If
\(r>0\), or if \(\chi(0,B,C)=0\), then
\begin{equation}
 \chi(A,B,A)=0 .
 \label{eq:three-reflected-boundary}
\end{equation}
Indeed, in the puncture basis substitute
\((\GA,\GB,\GD)=(U,V,V^{-1}U^{-1})\) in
\eqref{eq:betti-three-bd}. The reflected third loop becomes \(1\).
The boundary identity gives
\[
 u(U,V,V^{-1}U^{-1})=V^{-2s}f(V,U)^{-1}U^{-s}.
\]
Give \(U-1,V-1\) weights \(2,1\). The ordinary order of \(f-1\)
is \(n\) by \eqref{eq:weighted-magnus-cover} and the injectivity
of \(\iota\). Every word of \(f-1\) contains both letters, so
\(f(V,U)-1\) has weight greater than \(n\). In the covering basis
the pure \(\mathbf C\)-axis of \(u\) is \(\mathbf C^s\). It has
no coefficient through weight \(n\) under either hypothesis:
its possible weight-\(n\) term has \(B\)-count zero and survives
\(A=0\). If \(s\ne0\), this gives \(2^{1+v_2(s)}>n\),
also the initial weight of each scalar factor on the right after
subtracting \(1\). The factors are \(1\) when \(s=0\).
The substituted \(\mathbf C\) is \(V^{-1}U^{-1}V^{-1}\), with
initial term \(A+B^2\) after subtracting \(1\), where
\(A,B\) denote the initial forms of \(U-1,V-1\).
Reading weight \(n\) gives \(\widetilde\psi(A,B,A+B^2)=0\).
Taking \(B\)-count \(r\) gives \eqref{eq:three-reflected-boundary}.
Thus, for a leading class in the kernel on the larger graph, it
would suffice to derive independence of \(C\) from the defining
relations of \(\mathcal H_{123}\).
\end{remark}

\section*{Use of artificial intelligence}
My starting point was the double-cover interpretation and the
question of where the pentagon enters the filtration argument.
I used OpenAI's ChatGPT to draft and revise the manuscript and to
develop parts of the reconstruction and restricted tangent arguments
in Sections~\ref{sec:reconstruction}--\ref{sec:dd-comparison}.
Its mathematical contribution was
particularly substantial in the appendices. It helped derive and
write out the coefficient formulas from the Hirose--Sato deletion
operators in Appendix~\ref{app:coefficient-proof} and the
reconstruction of a leading component from its projection in
Appendix~\ref{app:dd-reconstruction}, using the deletion tests
and special identity of \cite{DD}.

\bibliographystyle{alpha}
\bibliography{inject}

\end{document}